\documentclass[12pt]{article}
\usepackage[latin1]{inputenc}
\usepackage{amsmath, amssymb, amstext, amsfonts, amsbsy, amsthm, fancyhdr, multicol,amscd}
\usepackage{graphicx}
\usepackage{amsfonts, amscd, amssymb}
\usepackage[all]{xy}
\usepackage[mathscr]{eucal}
\usepackage{amsthm,amsmath,amssymb,amscd} 
\usepackage{amsmath,amsthm,amssymb,graphicx}
\usepackage{epsfig}

\newtheorem{lemma}{Lemma}
\newtheorem{cor}{Corollary}
\newtheorem{Prop}{Proposition}

\newtheorem{Remark}{Remark}

\newcommand{\R}{\mathbb{R}}

\title{Dynamics in the Schwarzschild  isosceles three body problem}
\author{John A. Arredondo\footnote{Centro de Ciencias Matematicas, UNAM-Morelia, Mexico. Email: alexander@matmor.unam.mx}\,, Ernesto P\'erez-Chavela\footnote{Departamento de Matem\'aticas, Universidad Aut\'onoma Metropolitana-Iztapalapa,  Av. San Rafael Atlixco 186, Col. Vicentina, 09340 M\'exico, D.F., Mexico. Email: epc@xanum.uam.mx}\,,  Cristina Stoica\footnote{Department of Mathematics, 
Wilfrid Laurier University, Waterloo, Canada. Email: cstoica@wlu.ca}}

\begin{document}

\maketitle

\tableofcontents

\begin{abstract}
The  Schwarzschild potential, defined as 
$U(r)=-A/r-B/r^3$, where $r$ is the relative distance between two mass points and $A,B>0$, models astrophysical and stellar dynamics systems in a classical context.  In this paper we present a qualitative study of a three mass point system with mutual Schwarzschild interaction
 where the motion is restricted to isosceles configurations at all times. 
We retrieve the relative equilibria and provide the energy-momentum diagram. We further employ appropriate regularization transformations  to analyse the  behaviour of the flow  near triple collision. 

We emphasize the distinct features of the Schwarzschild model when compared to its Newtonian counterpart. We prove that, in contrast to the Newtonian case,  on any level of energy the  measure of the set on initial conditions leading to triple collision is positive. Further, whereas in the Newtonian problem triple collision is asymptotically reached only for zero angular momentum, in the Schwarzschild problem  the  triple collision is  possible for non-zero total angular momenta (e.g., when two of the mass points  spin infinitely many times around the centre of mass). This  phenomenon is known in celestial mechanics as the \textit{black-hole effect} and it is understood as  an  analogue in the classical context of the behaviour  near a Schwarzschild black hole.
 Also, while in the  Newtonian problem all  triple collision orbits  are necessarily homothetic,  in the Schwarzschild problem 
 this is not necessarily true. In fact,  in the Schwarzschild problem there exist  triple collision orbits which are neither homothetic, nor homographic.
 
\end{abstract}

\section{Introduction}

%
%
%
In 1916 Schwarzschild \cite{SCH} gave a solution of Einstein's field equations  which describes the gravitational field of a uncharged spherical non-rotating mass. 
%
%
It is known that the Schwarzschild metric leads - via a canonical formalism  that transposes the relativistic problem into the realm of celestial mechanics (see \cite{EDD}) - to a Binet-type equation, which describes the motion as governed by a force originating in a potential of the form 
\begin{equation}\label{sch.po}
U(r)=-\frac{A}{r}-\frac{B}{r^3}\,, \end{equation}
 where $r$ is the relative distance between two mass points and $A,B>0\,.$
%
%
The  potential above, which we call  the \textit{Schwarzschild potential}, was brought to the attention of the dynamics and  celestial mechanics communities  by Mioc et al. in \cite{MI, STM}. 
Classical dynamics in the Schwarzschild potential has interesting features, quite distinct when compared to their  Newtonian counterparts. In particular, collisions  may appear at non-zero angular momenta, giving rise to a so-called \textit{black-hole effect}, where a particle ``falls" into a Schwarzschild source field while spinning  infinitely many times around it. This is in contrast to the Newtonian $N$-body problem where  collision is possible only if the  total angular momentum is zero.  
%
%




 \smallskip
 The black hole effect was introduced in celestial mechanics by Diacu et al. in \cite{DMVS, D, DVS} (see also \cite{STT}).   For the  Schwarzschild  one-body  problem, the existence of the black hole effect was proven analytically  in \cite{STM} by employing a technique due to McGehee \cite{MG2} for the regularization of the vector field's singularity at collision. Since then, various studies concerning     Schwarzschild one- and two-body problems  have appeared (see \cite{MAB, MPS,  STT, VA}). 
Recently Campanelli et al. \cite{CA} 
have simulated numerically a three black hole system in the relativistic context, and observed that  
total collision is reached via spiral trajectories. One of the results of the present paper is proving  
that the same kind of dynamics is feasible in  a classical  mechanics model. 



\smallskip
We consider a particular case of the three body problem in which two equal point masses $m_1=m_2=M$  are confined to a horizontal  plane, symmetrically disposed with respect to their  common center of mass $O$, and  a third point mass $m$ is allowed to move only on the vertical  axis   perpendicular to the plane of masses $M$ through the point $O$.   
%
%
%
At any time, the configuration formed by the three mass points is that of an isosceles  triangle (possibly degenerated to a segment) and  the only rotations allowed are with respect to the vertical axis on which $m$ lies.  The dynamics is given, after taking into account the angular momentum conservation,  by a  two degree of freedom Hamiltonian system. It can be shown that for a three mass point system with two equal masses and with  rotationally invariant interactions, isosceles motions  form  a non-trivial  invariant manifold of the three mass point system phase space. We call   \textit{the isosceles Schwarzschild problem} the constrained three body problem as described above where the mutual interaction between the mass points is given by a Schwarzschild potential of the form (\ref{sch.po}).

%

\smallskip
 Isosceles three body problems are often considered as case studies  for the more complicated  dynamics of the $N$-body problem. 
  One of the main references is a study by Devaney \cite{DE} where the author employs  McGehee's technique \cite{MG}  and presents a qualitative description of the planar isosceles three body problem (i.e., the isosceles three body problem with zero angular momentum) with emphasis on orbits which begin or end in a triple collision. In \cite{MO} Moeckel  uses geometrical methods to construct an invariant set containing a variety of periodic orbits which exhibit close approaches to triple collision and wild changes of configuration. He also finds heteroclinic connections between these periodic orbits, as well as oscillation  and capture orbits.  Simo and Martinez \cite{SM} apply analytical and numerical tools to study  homoclinic and heteroclinic orbits which connect triple collisions to infinity, and use  homothetic solutions to obtain a characterization of the orbits which pass near triple collision.  In \cite{EB}  ElBialy uses  (rather than the kinetic mass matrix norm)   the Euclidian norm to perform a  McGehee-type change of  coordinates in order to discuss the flow behavior near the collision manifold as the ratio between the   $m/M \to 0.$ 
  In a recent study, Mitsuru and Kazuyuri \cite{MK}  deduce numerically the existence of infinite families of relative periodic orbits.

 
 

\smallskip
The first part of our study follows classical methodology: we write the Hamiltonian, use the rotational symmetry to reduce the system, and obtain  
the reduced Hamiltonian as a sum of the kinetic energy and the reduced (amended) potential. 
The internal parameters are the energy $h$ and the total angular momentum $C$. We retrieve the relative equilibria, that is,  solutions where the three mass points  are steadily rotating about their common centre of mass, and discuss their stability modulo rotations. 

 \smallskip
Recall that in the Newtonian three body problem  there are two classes of relative equilibria (up to  a homothety): the Lagrangian  ones, where the three mass points   form an equilateral triangle and the rotation axis is perpendicular to the plane determined by the mass points; and the Eulerian ones, where the mass points  are in a collinear configuration and the axis of rotation is perpendicular to the line of the mass points at the centre of mass.  For certain mass ratios the Lagrangian relative equilibria are linearly stable, whereas the Eulerian relative equilibria are always unstable. 
In the  isosceles problem, since the axis of rotation is  perpendicular to the line of the equal masses at their centre of mass, one retrieves only the Eulerian relative equilibria, which are still unstable, even when  considered only in the invariant manifold  of  isosceles motions. In the Schwarzschild isosceles problem we find  no relative equilibria  for small angular momenta. For  angular momenta  larger than a critical value,  we find two distinct Eulerian relative equilibria: one unstable, similar to the Newtonian case, and one stable modulo rotations about the vertical axis. This is the case even if  $B$ is small, that is when   Schwarzschild may be considered as  a perturbation of the Newtonian model. 
 
 \smallskip
 The second part of our study focusses on motions near collisions  when $m<<M.$ A similar analysis may be performed for general ratios $m/M$, which we expect would lead to similar conclusions, but we leave this   for a future study. By  using a  McGehee-type transformation (similar as in \cite{DE}) we introduce coordinates which  
regularize   the flow at  triple  collisions.  The triple collision singularity  appears as a  manifold  pasted into the phase space for all levels of energy and angular momenta. This manifold contains fictitious dynamics which is  used to draw conclusions about the motions near, and/or leading to, triple collisions.
In our case, the triple collision  manifold is not closed and has two edges which correspond to triple collisions attained through motions in which  the equal masses are a black-hole type binary collision and the third  mass  $m$ is on one side of the vertical axis (see Figure \ref{orbitas13}). The condition $m<<M$ insures the presence of six equilibria on the triple
collision manifold, three corresponding to orbits which begin on the manifold, called \textit{ejection}  orbits, and three which end on it, called \textit{collision} orbits.   We  analyse the flow on the triple collision manifold, including the  stability of equilibria, and we  deduce that    on every energy level the set of initial conditions ejecting from or leading to triple collision has positive Lebesgue measure.  We also find a condition over the set of parameters so that certain orbit connections are satisfied. 

\smallskip
It is important to remark that one of the main differences between the Newtonian and Schwarzschild isosceles problem at triple collisions is that in the latter the \textit{binary collisions do not regularize as elastic bounces}, but as  one dimensional invariant manifolds. This is not surprising, taking into account previous results on the two body  collision in the  problem with non-gravitational interactions (see \cite{MG2}; also \cite{STT}). 

 \smallskip
It would be interesting to adapt the McGehee-type transformation used by ElBialy  \cite{EB}  to   the context of Schwarzschild problem. We reckon that this would allow the analysis of the flow at collision for both $m/M \neq 0$ and $m/M=0$ (i.e., $m=0$) and  that, similarly to in the Newtonian case, the collision manifold will have different topologies for non-zero and zero mass ratios.   However, taking into account the remark in the paragraph above, it is hard to predict if other similarities would occur when about the discussion of the near-collision flow.  We defer this problem to a future investigation.

 \smallskip
We continue by studying homographic solutions and, implicitly,  central configurations. By definition,  a \textit{homographic} solution for  a $N$ mass-point system is a solution along which the  geometric configuration of the mass points is similar to the  initial geometric configuration.
There are two extreme cases of homographic motion: if the motion of the mass points is  a steady rotation, then the solution is in fact a  relative equilibrium; and if the mass points evolve on straight lines through the common centre of mass, then the solution is called \textit{homothetic}.  
The geometric configuration (within its similarity class) of the points  of a homographic solution is called a  \textit{central configuration} if at all times the position vectors are parallel to the acceleration vectors.
 It can be shown that, for rotationally-invariant potentials, there are only two instances where central configurations are possible: either the mass points are in a relative equilibrium, or they are homothetic. Central configurations play an important role in understanding $N$-body systems. Excellent references for this subject can be found in   \cite{Saa} and \cite{Saa05}.  For the Newtonian problem, the mass points tend to such central configurations as they approach total collisions; it is worth to mention that an outstanding open problem  is the finiteness of central configurations in the Newtonian  $N$-body problem (see \cite{Smale}). 

 \smallskip
We show that for the Schwarzschild isosceles problem,  homographic  motions are  confined  to the horizontal plane and with $m$ resting at $O$ for  all times. In particular, the only central configurations are given by the Eulerian relative equilibria and the (collinear) homothetic motions.  This is in agreement with the analysis of   central configurations for the   three-body problem with generalized Schwarzschild interaction presented by Arredondo et al. in \cite{AP}. 
\begin{table}[h!]
  \centering 
  \begin{equation*}
\begin{array}{|c|c|}
\hline \\ \text{The Newtonian  isosceles problem} & \text{The Schwarzschild   isosceles problem}\\
\,\\
 \hline \\
   \,
& \text{There are two  collinear (Eulerian) relative}
\\
 \text{There is one  collinear (Eulerian) relative}   &  \text{equilibria (up to a permutation  of the}
\\
\text{ equilibrium (up to a permutation  of the }   & \text{equal masses), one  stable (modulo rotations)}, \\
\text{equal masses) which is unstable.} &  \text{ and one unstable.}\\
\,\\
 \hline \\
 \text{On every level of energy, the set of} &  \text{On every level of energy, the set of}  \\
 \text{initial conditions  leading to triple collision}   &  \text{initial conditions leading to triple collision}\\
  \text{has zero  Lebesgue measure.}   &  \text{has positive  Lebesgue measure.}
 \\
 \,\\
 \hline \\
\,& \text{Triple collision is possible for all $C$.} \\
  \text{Triple collision is possible only when $C=0$.}  & \text{For $C\neq 0$, the equal masses display}\\
 \,& \text{black-hole type motion.}
 \\
 \hline \\
 \text{All triple collision orbits are homothetic.} & \text{There are triple collision orbits} \\
\text{\,} & \text{which are not homothetic.}
 \\
 \text{\,} & \text{(Also, there are triple collision orbits}
 \\
 \text{\,} & \text{which are not homographic.)}
 \\
 \,\\
 \hline\\
\, & \text{There are asymptotic geometric configurations} \\
 \text{All asymptotic geometric configurations}  & \text{at triple collision  which are not}\\
\text{at collision are central configurations.} & \text{central configurations. This is true}\\
\, & \text{for both $C=0$ and $C\neq 0$ cases.}\\
 \,\\
  \hline
\end{array}
\end{equation*}
\caption{Newtonian versus Schwarzschild dynamics in the isosceles problem.}
  \label{table}
\end{table}

 \smallskip
Finally we analyse  motions near triple collisions.  We discuss the asymptotic (limiting) geometric configurations of the solutions corresponding to the ejection/collision orbits as they depart from or tend to triple collision. In the Newtonian problem, such limiting configurations are  associated to a central configurations as the solutions corresponding to the ejection/collision orbits  are homothetic. We prove that this is not the case in the generic Schwarzschild problem. (The non-generic case is given by  condition (\ref{generic_pot}). See also  Section \ref{hom_sol}.) Moreover, there are limiting triangular  configurations which are not even associated to homographic solutions. To our knowledge, this is the first time when such ``non-homographic"  configurations are observed to be  limiting configurations at triple collision.  We further remark that on any energy level, the set of initial conditions leading to triple collision with a  limiting  geometric configuration of a homographic  solution is of positive Lebesgue measure. These solutions correspond to  motions starting/ending in triple collision where mass $m$ on the vertical axis crosses the horizontal plane  infinitely many times before collision (i.e., $m$ oscillates about the centre of mass of the binary equal mass system). 
The set of initial conditions leading to triple collision with a non-degenerate  non-homographic triangular limiting configuration  is of zero  Lebesgue measure.  
In all cases, whenever the angular momentum is not zero, collisions are attained while the equal masses perform black-hole type motions. We end by proving that for negative energies, there is an open set of initial conditions for which solutions end in double (i.e., the collision of the equal masses while $m$ is above or below the horizontal plane) or in triple  collisions with $m$  crossing the horizontal plane a finite number of times.




\medskip
Our study emphasizes the differences between the Newtonian and the  Schwarzschild model. Our findings are summarized in Table \ref{table}.  As mentioned, one of the most important dissimilarities  concerns the approach to  total collapse. Besides the presence of the black hole effect, the Schwarzschild problem displays asymptotic total collision trajectories which are not homographic; in particular, this implies that there exist asymptotic geometric configurations at triple collision which are not central configurations.  

\medskip
To put these conclusions in context, we note that given the presence of the ``strong-force'' $-B/r^3$ term in the Schwarzschild potential which dominates at small distances,  leads to the expectation that  black hole effects near collisions are present (see also \cite{STT}); 
 initially, the main goal of this study was to prove this for three mass point interactions.  The existence of the non-homographic triple collisions orbits is  due to the non-homogeneity of the potential, as it is given by a sum of two homogeneous terms which are taken in a ``generic" position (see Section \ref{V_and_W}). 
 %
 Related work  was performed on the three-body problem with quasi-homogeneous interaction, a generalization of the Schwarzschild potential of the form $-A/r^a- B/r^b,$ $1\leq a<b$.  Diacu \cite{DF}  found that in the quasi-homogeneous three body problem,  the set of collision orbits form asymptotically \textit{quasi-central configurations}, that is, geometric configurations of orbits which are homographic only with respect to  the term $- B/r^b$ of the potential.  In \cite{E}, Diacu et al.  studied  the so-called \textit{simultaneous central configuration} for quasi-homogeneous interactions; these are  central configurations arising  in the non-generic case (see equation  (\ref{generic_pot}))  which we do not consider here.    Perez-Chavela et al. \cite{PV} studied the collinear quasi-homogeneous three body problem (it is assumed that there is no rotation) and proved that the set of initial conditions leading to triple collision has positive Lebesgue measure.  In the same paper the authors show that there are triple collision orbits which are not asymptotic to a central configuration, but these orbits do not originate/terminate in a fixed point on the collision manifold. (They involve infinitely many double collisions of a pair of outer masses.)
%



\medskip
The paper is organized as follows. In  Section \ref{SIP.sect} we describe the Schwarzschild isosceles problem and reduce the system to a two degrees of freedom Hamiltonian system. We further study the relative equilibria and their stability, and provide the energy-momentum bifurcation diagram. In  Section \ref{CM.sect} we introduce new coordinates to regularize  singularities due to double and triple collisions, define the triple collision manifold and describe the flow  behavior on it. In  Section \ref{GFB.sect} we study homographic solutions. Finally, in Section \ref{n_c_o} we analyse the  triple collision/ejection orbits.


\section{The Schwarzschild isosceles problem}\label{SIP.sect}

Consider three point masses with masses $m_1=m_2=M$ and $m_3=m$ interacting  mutually via a  Schwarzschild-type potential.  Let  $ {\bf r}_1$ and ${\bf r}_2$  the positions vectors in Jacobi coordinates, that is  $ {\bf r}_1$ is the vector from the particle of mass $m_1$ to the particle of mass $m_2$ and ${\bf r}_2$  is the vector from the center of mass of the first two particles to the particle with mass $m_3$. The associated momenta are ${\bf p}_1$ and ${\bf p}_2$. In these coordinates the respective Hamiltonian is given by
   \begin{align}\label{H1}
  H= \frac{1}{m} & {\bf p}^2_1 + \frac{2M+m} {4Mm} {\bf p}^2_2 +
  U_{12}( \left | {\bf r}_1 \right|  ) + U_{13}\left( \left | {\bf r}_2 + \frac{1}{2} {\bf r}_1 \right| \right) +
   U_{23}\left( \left |  {\bf r}_2 - \frac{1}{2} {\bf r}_1\right | \right),
  \end{align}
where the  $U_{ij}$, $i,j=1,2,3$, $i\neq j,$ are  Schwarzschild type  potentials. The system  is invariant under the diagonal action  of the $SO(3)$ group of spatial rotations  on the configuration space $ \mathbb{R}^3\times \mathbb{R}^3 \setminus \{ ({\bf r}_1 \,,{\bf r}_2) \,|\, {\bf r}_1\neq 0 \}$ 
which leads to the conservation of the level sets of the momentum map
\[J({\bf r}_1, {\bf r}_2, {\bf p}_1, {\bf p}_2) =  {\bf r}_1 \times {\bf p}_1 + {\bf r}_2 \times {\bf p}_2.\]
%
%


In our modeling, the equal point masses $M$  are  confined to a horizontal  plane and are symmetrically disposed with respect to their  common center of mass $O$, and  $m$ is allowed to move on the vertical  axis perpendicular to the $xy$ plane in $O.$ For motions with zero angular momentum, the three masses lie in their initial plane for all times, whereas for   motions with non-zero angular momentum  the masses $M$   are rotating about the vertical axis on which $m$ lies.
The motion is described by a  Hamiltonian system which in coordinates  ${\bf r}_1=(x_1,y_1,0)$, ${\bf r}_2=(0,0,z_2)$, and momenta ${\bf p}_1=(p_{x_1}, p_{x_2}, 0)$, ${\bf p}_2=(0,0, p_{z_2})$, respectively, is
     \begin{align}
     &H: \left(\mathbb{R}^3\setminus \{(x_1,y_1,z_2)\,|\, x_1^2+y_1^2=0\}\right)\times \mathbb{R}^3 \to \mathbb{R} \nonumber\\
  &H(x_1,y_1,z_2, p_{x_1}, p_{x_2},  p_{z_2})= \frac{1}{M} \left(p^2_{x_1}+p^2_{y_1} \right) + \frac{2M+m} {4Mm} p_{z_2}^2 +
  U(x_1, y_1, z_2 ), \label{H2}
  \end{align}
where the potential has the form
       \begin{equation}\label{pot2}
   U(x_1, y_1, z_2 )= -\frac{A}{\sqrt{x_1^2+y_1^2}}- \frac{B}{ \sqrt{(x_1^2+y_1^2)^3}} - \frac{4A_{1}}{\sqrt{x_1^2+y_1^2 +4z_2^2}}-\frac{16B_1}{\sqrt{(x_1^2+y_1^2+4z_2^2)^3}}\,.
    \end{equation}
\noindent
The angular momentum  integral is given by
\begin{equation}
C= x_1 p_{x_2}-x_2 p_{x_1}.
\end{equation}

\begin{Remark}
In the isosceles Schwarzschild problem there are six (external) parameters: $M$, $m$, and $A$, $A_1,$ $B$ and $B_1.$   Since without loosing generality, we could take one of parameters  to be one (e.g., one of the masses), there are  five independent parameters.
 
 \end{Remark}

\noindent
It is convenient to pass to cylindrical coordinates $(x_1, y_1, z_2, p_{x_1}, p_{y_1}, p_{z_2}) \to (R, \phi, z, P_R, P_{\phi}, P_z)$ given by the change of coordinates
\[x_1= R \cos \phi, \quad\ y_1= R \sin \phi, \quad z_2=z,\]
and its associated (canonical) transformation of the momenta.
The Hamiltonian becomes
 \begin{equation}\label{H3}
  H(R, \phi, z, P_{R}, P_{\phi}, P_z)= \frac{1}{M} \left( P_R^2+ \frac{P_{\phi}^2}{R^2}\right) + \frac{2M+m} {4Mm} P_{z}^2 +
  U(R, z ),
  \end{equation}
with
 \begin{equation}\label{P_cyl}
U(R, z)= -\frac{A}{R}-\frac{B}{R^3}-\frac{4A_{1}}{\sqrt{R^2+4z^2}}-\frac{16B_1}{(R^2+4z^2)\sqrt{R^2+4z^2}}.
\end{equation}
The equations of motion for the variables $(\phi, P_{\phi})$ are

\[ \dot \phi =\frac{\partial H} {\partial P_{\phi}}=   \frac{2 P_{\phi}}{MR^2}, \quad \quad \dot p_{\phi} = - \frac{\partial H} {\partial \phi}=0.\]

\noindent
leading to the explicit equation of the angular momentum conservation  
\[P_{\phi}(t)=const.=:C\,.\] 
Using the above,  we obtain  a two degree of freedom Hamiltonian system  determined by  the reduced Hamiltonian
  \begin{equation}\label{H_red_1}
  H_{\text{red}}(R,  z, P_{R},  P_z; C):= \frac{1}{M} \left( p_R^2+ \frac{C^2}{R^2}\right) + \frac{2M+m} {4Mm} P_{z}^2 +
  U(R, z )\,,
  \end{equation}
that is, a system of the form ``kinetic + potential":
 \begin{equation}\label{H_red}
  H_{\text{red}}(R,  z, P_{R},  P_z; C)
=\frac{1}{2}(p_R  \quad  p_{z}) \left(
\begin{array}{cc}
\frac{2}{M} &0\\
0&\frac{2M+m} {2Mm}
\end{array}
\right)
\left(
 \begin{array}{c}
P_R\\
P_z
\end{array}
 \right)
+
  U_{\text{eff}}(R, z ),
  \end{equation}
with the effective (or amended) potential given by
   \begin{equation}\label{U_eff}
   U_{\text{eff}}(R, z; C):=\frac{C^2}{MR^2}  + U(R,z) = \frac{C^2}{MR^2}  -\frac{A}{R}-\frac{B}{R^3}-\frac{4A_{1}}{(R^2+4z^2)^{1/2}}-\frac{16B_1}{(R^2+4z^2)^{3/2}}\,.
   \end{equation}
The equations of motion are
   \begin{align*}
   \dot R &= \frac{2}{M}P_R\,, \quad \quad \quad \quad
   \dot P_r= - \frac{\partial H}{\partial R} = - \left( - \frac{2C^2}{MR^3} +\frac{A}{R^2}  + \frac{3B}{R^4}  + \frac{4 A_1R }{(R^2+4z^2)^{3/2}}
 +\frac{3 \cdot 16B_1 R}{(R^2+4z^2)^{5/2}}
   \right),\\
  \dot z&= \frac{2M+m}{2Mm}P_z\,,  \quad \,\,\dot P_z= -  \frac{\partial H}{\partial z}= - \left(   \frac{4 A_1}{(R^2+4z^2)^{3/2}}  +    \frac{3 \cdot 16 B_1 }{(R^2+4z^2)^{5/2}}  \right)4z\,.
   \end{align*}
   Along any integral solution, the energy is conserved:
   \begin{equation}\label{H_red_eff}
H_{\text{red}}\left(R(t),  z(t), P_{R}(t),  P_z(t); C\right)=const.=h.
\end{equation}

 \begin{Remark} \label{planar_motions}
 The submanifold 
\begin{equation}
\label{inv_plane}
 \{(R, z, P_R, P_z) \in (0, \infty) \times \mathbb{R} \times \mathbb R  \times  \mathbb{R} \,|\, z=0, P_z=0\}
 \end{equation} 
is invariant. Physically, this submanifold contains planar motions, with the two masses $M$ symmetrically disposed with respect to their midpoint $O$ where $m$ rests at all times. The motion on this manifold is the subject of Section \ref{coll-eject}.

 \end{Remark}


\subsection{Relative equilibria}\label{Sect_RE}

Following classical methodology, for non-zero angular momenta, the equilibria of (\ref{H_red}) are in fact relative equilibria, that is dynamical solutions  which are also one-parameter orbits of the symmetry group.  In our case, relative equilibria correspond to trajectories  where the mass points $M$ are steadily rotating about the vertical $z$ axis. Note that since $m$ lies on the $z$ axis, it does not  ``feel"  such rotations. 


\begin{Prop}
\label{real_equilibria}
Consider  the spatial isosceles Schwarszchild three body problem   and let  $C$ be the magnitude of the angular momentum. Without loosing generality let us consider $C>0.$ (For $C<0$ the same results are obtained, but where spin of the angular momentum  is reversed.) Denote 
\begin{equation}
\label{notation-roots}
\alpha:=M(A+4A_1)\,,\quad \quad\beta:=M(B+16B_1)\,.
\end{equation}
and let
\begin{equation}
\label{c_zero_not}
C_0:=\sqrt[4]{3\alpha \beta}\,.
\end{equation}
 Then:

\begin{enumerate}
\item If $C< C_0$ then there are no relative equilibria.

\item If $C=C_0$ then there is one relative equilibrium, and it is of collinear configuration
with the equal point masses situated at
\begin{equation}\label{col_eq}
R_{0}=\frac{C^2}{\alpha}.
\end{equation}
This relative equilibrium is of degenerate stability, having  a zero pair of eigenvalues.

\item If $C> C_0$ then there are two relative equilibria, both of collinear configuration with the equal point masses situated at $(R, z)=(R_i,0)$, $i=1,2$ where
\begin{equation}\label{e23}
R_{1}=\frac{C^2+\sqrt{C^4-C_0^4}}{\alpha} \quad  and \quad R_{2}=\frac{C^2-\sqrt{C^4-C_0^4}}{\alpha}.
\end{equation}
The relative equilibrium $(R_{1}, 0)$  is non-linearly stable modulo rotations, whereas $(R_{2}, 0)$ is unstable.
\end{enumerate}
\end{Prop}

\noindent
\textbf{Proof:} The relative equilibria of the system given by the Hamiltonian (\ref{H_red}) correspond to  the critical points of the effective potential (\ref{U_eff}):
\begin{equation}\label{e22}
\frac{\partial U_{\text{eff}}}{\partial R} =\frac{\partial U_{\text{eff}}}{\partial z} =0.
\end{equation}

For $C\geq C_0=\sqrt[4]{3 \alpha \beta}$ we find  solutions with $z=0$ and $R=R_{1, 2}$, given by the roots of 
\begin{equation}
\label{eq-roots}
\alpha R^2 -2C^2 R +3\beta=0\,.
\end{equation}
We have
\begin{equation}\label{e_prim_23}
R_{1,2}= \frac{C^2\pm\sqrt{C^4-3 \alpha \beta}}{\alpha}= \frac{C^2\pm\sqrt{C^4-C_0^4}}{\alpha}.
\end{equation}
For $C>C_0,$ the nonlinear stability  of the relative equilibria  may be established by calculating  $D_2U_{\text{eff}}$  at $(R_{1, 2}, 0)$.  More precisely, if $D_2U_{\text{eff}}$ at one of the relative equilibrium is positive definite, then the particular relative equilibrium is nonlinear stable modulo rotations (see \cite{A, ME}).
We have
\[D_2U_{\text{eff}}|_{\text{z=0}}=
\left(
  \begin{array}{cc}
    \dfrac{6C^2}{MR^4}-\dfrac{2(A+4A_1)}{R^3}-\dfrac{12(B+16B_1)}{R^5} & 0 \\
    0 & \dfrac{16A_1}{R^3}+\dfrac{192B_1}{R^5} \\
  \end{array}
\right)\,.
\]
Now we check  the positive definiteness of $D_2U_{\text{eff}}$  at $(R_{1}, 0)$ and $(R_{2}, 0)$. For this we have to analyse the behavior of the first entry in the matrix above. Let us define
\begin{equation}\label{e24}
f(R)=-2\alpha R^2+6C^2R-12\beta,
\end{equation}
and note that
\[ \dfrac{6C^2}{MR^4}-\dfrac{2(A+4A_1)}{R^3}-\dfrac{12(B+16B_1)}{R^5} = \frac{f(R)}{MR^5}\,.
\]
From equation (\ref{eq-roots}) we have that any of the roots $R=R_{1, 2}$ verify
 \[R^2=\dfrac{2C^2R-3\beta}{\alpha}.\]
Substituting the above into $f(R)$, after some calculations we have
 %
%
%
%
\begin{equation}\label{e_prim_24}
f(R_{1, 2})=\frac{2C^2(C^2\pm\sqrt{C^4-3\alpha\beta})-6\alpha \beta }{\alpha}.
\end{equation}
At $R_{1}$,   we have
\begin{align*}
f(R_{1})&= 
\frac{2(C^4-3\alpha \beta)}{\alpha} + \frac{2C^2\sqrt{C^4-3\alpha \beta}}{\alpha}>0
%
\end{align*}
and so $D_2U_{\text{eff}}|_{z=0,R=R_{1}}$ is definite positive, i.e.,  the relative equilibrium $(R_{1},0)$ is nonlinear stable.

At $(R_{2},0)$, the Hessian matrix is indefinite, and so it does not give information about its stability. We calculate then the  spectral stability of $(R_{2},0)$  by computing
the eigenvalues of the  matrix linearization:
\[L=
\left. \left(
  \begin{array}{cccc}
    0 & 0 & \frac{2}{M} & 0 \\
    0 & 0 & 0 & \frac{2M+m}{2Mm} \\
    -\frac{\partial^2 U_{\text{eff}}}{\partial R^2}  & -\frac{\partial^2 U_{\text{eff}}}{\partial R\partial z} & 0 & 0 \\
    -\frac{\partial^2 U_{\text{eff}}}{\partial R\partial z} & -\frac{\partial^2 U_{\text{eff}}}{\partial z^2}  & 0 & 0 \\
  \end{array}
\right) \right|_{(R_{2}, 0)}=
\left(
  \begin{array}{cccc}
    0 & 0 & \frac{2}{M} & 0 \\
    0 & 0 & 0 & \frac{2M+m}{2Mm} \\
    -\frac{f(R_2)}{MR_2^5}  & 0 & 0 & 0 \\
    0 & -\left( \dfrac{16A_1}{R_2^3}+\dfrac{192B_1}{R_2^5} \right)  & 0 & 0 \\
  \end{array}
\right)\,.
\]
At $z=0$, these correspond to:
\[\lambda_{1,2}=\pm 4i\sqrt{\frac{2M+m}{2Mm}\left( \dfrac{16A_1}{R_2^3}+\dfrac{192B_1}{R_2^5} \right) }\]
and
\[
\lambda_{3,4}=\pm \sqrt{\frac{2}{M}  \left(   -\frac{f(R_2)}{MR_2^5} \right)} = \pm\frac{2}{M} \sqrt{
\frac{1}{R_2^5} (3 \beta - C^2 R_2)
}.
\]
The eigenvalues $\lambda_{1,2}$ are purely imaginary. A direct calculation shows that $3 \beta - C^2 R_2>0$ and so  $\lambda_{3,4}$ are real. In conclusion,  the relative equilibrium $(R_2, 0)$ is unstable. $\square$


\subsection{Energy-momentum diagram}\label{Energy}

 The energy-momentum diagram provides the location of the relative equilibria  in the $(h,C)$ parameter space.   As known (see \cite{SSM}), this is the set of points  where the topology of the phase space changes.

\begin{figure}[h]
\centerline
{\includegraphics[scale=1]{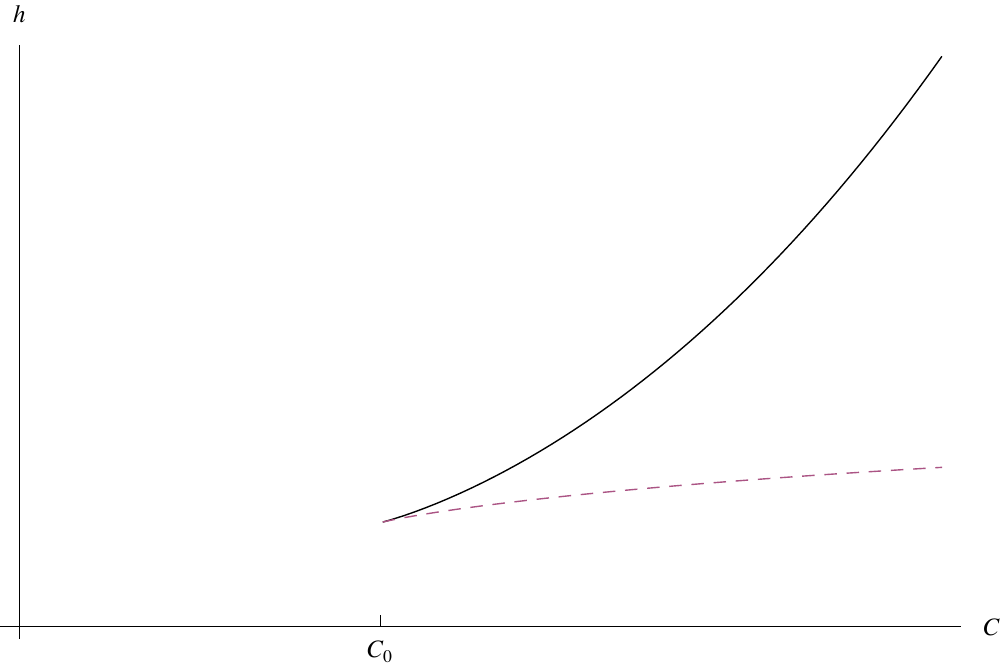}}
{\caption{ \label{energy_mom_fig}The generic graph of the energy momentum map (where $C>0$).  The continuous line corresponds to  the stable relative equilibria  $(R_1, 0)$, whereas the  dashed  line corresponds to the unstable relative equilibria $(R_2, 0).$ 
(The figure is generated for $M=1$,  $A=A_1=1$, $B=B_1=0.2$, that is $\alpha=5$ and $\beta=3.4$) 
}}
\end{figure}

 In our case, the energy-momentum curve is determined   by eliminating $R_{i},$ $i=1,2,$ as given by formula (\ref{e23}) from  the energy relation
  at a relative equilibrium (where $P_R=P_{z}=0$ and $z=0$):
  \begin{align}\label{energy_RE}
\frac{C^2}{MR_{i}^2}  -\frac{A+4A_1}{R_{i}}-\frac{B+16B_1}{R_{i}^3}=h.
  \end{align}
or, using the notation (\ref{notation-roots})
\begin{align}
\frac{C^2}{R_{i}^2}  -\frac{\alpha}{R_{i}}-\frac{\beta}{R_{i}^3}=h.
  \end{align}%
On this curve there are  two points where the relative equilibria curves intersect. The momentum $C$ of these points are given by the equation
\[
\frac{C^2}{R_{1}^2}  -\frac{\alpha}{R_{1}}-\frac{\beta}{R_{1}^3} = \frac{C^2}{R_{2}^2}  -\frac{\alpha}{R_{2}}-\frac{\beta}{R_{2}^3}\,,
\]
which we can re-write as 
\[
   \left(\frac{1}{R_1} -  \frac{1}{R_2} \right) \left[ C^2 \left( \frac{1}{R_1} +  \frac{1}{R_2}   \right)  - \alpha 
-\beta\left( \frac{1}{R_1^2} +  \frac{1}{R_1R_2}+ \frac{1}{R_2^2}   \right) \right]=0\,.
\]
An immediate solution is $C=C_0$ where $R_1=R_2=R_0.$ A second solution  is given by  the equation 
\[C^2 \left( \frac{1}{R_1} +  \frac{1}{R_2}   \right)  - \alpha 
-\beta\left( \frac{1}{R_1^2} +  \frac{1}{R_1R_2}+ \frac{1}{R_2^2}   \right)=0\]
which, after substituting the formulae (\ref{e23}) for $R_{1,2}$, leads to $C= \sqrt[4]{3 \alpha\beta}=C_0.$ So we deduce that the energy-momentum curve has no self-intersections no-matter the choice of the parameters.  
A generic graph of the energy momentum map  is presented in Figure (\ref{energy_mom_fig}).

\begin{Remark}
Recall that in the Newtonian case (i.e., when $B=B_1=0$) the isosceles problem displays   collinear  (Eulerian) relative equilibria, which are unstable.  We observe that in the presence of   the inverse cubic terms there are two families of collinear relative equilibria and that one of these is nonlinearly stable. This is the case even if $B$ and $B_1$ are small, and so the inverse cubic terms can be thought of as a perturbation of the Newtonian problem. 

\end{Remark}


\section{The triple collision manifold}\label{CM.sect}
In this section we regularize  the equations of motion of the isosceles  Schwarzschild three body problem so that the dynamics at triple and double-collisions appear on a  fictitious \textit{collision} invariant manifold. We further discuss the orbit behaviour on the collision manifold.
From now on, unless otherwise stated, we assume that $M>>m.$ 

\subsection{New coordinates}

To start the study of the dynamics near singularities  (i.e., collisions) it is convenient to transform the system associated to the Hamiltonian (\ref{H_red}), so that the singularities are regularized.  For this we follow closely the McGehee technique as used in the  Newtonian isosceles  problem by Devaney (see \cite{DE}). Denoting
\[{\bf x}:=
\left(
\begin{array}{c}
R\\
z
\end{array}
\right)
, \quad  {\bf p}:=
 \left(
\begin{array}{c}
P_R\\
P_z
\end{array}
\right), \quad \text{and}\,\,\,\,T=
\left(
\begin{array}{cc}
\frac{M}{2} &0\\
0&\frac{2Mm}{2M+m}
\end{array}
\right)\,,
\]
we introduce the  coordinates $(r,v, {\bf {s}}, {\bf u} )$  defined by
\begin{align}\label{sys0}
r&=\sqrt{{\bf x}^{t}T{\bf x}},\quad \quad \quad
v=r^{\frac{3}{2}}({\bf s}\cdot {\bf p}),\\
{\bf s}&=\frac{{\bf x}}{r},\quad \quad \quad  \quad \quad \,\,\,
{\bf u}=r^{\frac{3}{2}}(T^{-1}{\bf p}-({\bf s}\cdot {\bf p}) {\bf s}).\nonumber
\end{align}
%
%
%
Note that $r=0$ corresponds to $R=z=0$, i.e., to the triple collision of the bodies. One may verify that in the new coordinates we have
%
%
%
%
that ${\bf s}^tT{\bf s}=1$ and ${\bf s}^tT{\bf u}=0$. The equations of motion read

\begin{align*}
\dot{r}&=\frac{v}{r^{\frac{3}{2}}},\\
\dot{v}&=\frac{3}{2}\frac{v^2}{r^\frac{5}{2}}+\frac{{\bf u}^t T {\bf u}}{r^\frac{5}{2}}+\frac{1}{r^\frac{3}{2}}\frac{2C^2}{Ms_{1}^2}-\frac{1}{r^\frac{1}{2}}V({\bf s})-\frac{3}{r^\frac{5}{2}}W({\bf s}),\\
\dot {\bf s}&=\frac{{\bf u}}{r^{\frac{5}{2}}},\\
\dot {\bf u}&= \frac{1}{2}\frac{v}{r^\frac{5}{2}} {\bf u}+ \left(-\frac{{\bf u}^t T {\bf u}}{r^\frac{5}{2}}-\frac{2C^2}{Ms_{1}^2r^\frac{3}{2}}+\frac{1}{r^\frac{1}{2}}V({\bf s})+\frac{3}{r^\frac{5}{2}}W({\bf s})\right) {\bf s}\\
\,\\
& + \frac{1}{r^\frac{1}{2}} \left(
                            \begin{array}{c}
                              \dfrac{2}{M}\dfrac{\partial V}{\partial s_1}  \\
                                         \dfrac{2M+m}{2Mm}\dfrac{\partial V}{\partial s_2}  \\
                                       \end{array}
                                     \right)+
                                   \frac{1}{r^\frac{3}{2}} \left(
                                     \begin{array}{c}
                                     \dfrac{\partial }{\partial s_1} \left(  - \frac{2 C^2}{M^2s_1^2} \right)\\
                                     0
                                       \end{array}
                                     \right)+
                                    \frac{1}{r^\frac{5}{2}} \left(
                                       \begin{array}{c}
                                         \dfrac{2}{M}\dfrac{\partial W}{\partial s_1}  \\
                                         \dfrac{2M+m}{2Mm}\dfrac{\partial W}{\partial s_2}  \\
                                       \end{array}
                                     \right),
\end{align*}
where
\[V({\bf s} ) =\frac{A}{s_1}+\frac{4A_1} {\left( s_1^2+4s_2^2 \right)^{1/2} } \quad \quad\text{and} \quad \quad
W({\bf s}) = \frac{B}{s_1^3}+\frac{16B_1}{\left( s_1^2+4s_2^2 \right)^{3/2}}.\]

\medskip
\noindent
We further introduce the change of coordinates given by

\[{ \bf s}=\sqrt{(T^{-1})}(\cos\theta,\sin\theta)^t  \quad \text{and} \quad {\bf u}=u\sqrt{(T^{-1})}(-\sin\theta,\cos\theta)^t\]

\medskip
\noindent
where $\displaystyle{-\dfrac{\pi}{2}<\theta<\dfrac{\pi}{2}}$ so that the boundaries $\displaystyle{\theta=\pm \frac{\pi}{2}}$ correspond in the original coordinates to $R=0,$ that is, to  double collisions of the masses $M.$ More precisely, at $\theta=\pi/2$ we have $R=0$ and $z>0$, whereas at $\theta=-\pi/2$, $R=0$ and $z<0$. One may easily verify that  ${\bf u}^t T {\bf u}= u^2$ and 
$\displaystyle{
\dot {\bf u}=({\dot u}/{u})   {\bf u}  - u \,\dot \theta \,{\bf s}.}
$
%
Denoting
 \begin{equation}\label{mu}
 \mu:=\dfrac{2M+m}{m}\,
 \end{equation}
and applying the time re-parametrization  $dt = r^{\frac{5}{2}}d\tau$, we obtain the system
\begin{align}\label{sys2}
r^{\prime}&=rv,\nonumber \\
v^{\prime}&=\frac{3}{2}v^2+u^{2}+\frac{C^2}{\cos^2\theta}r-r^2V(\theta)-3W(\theta), \\
\theta^{\prime}&=u ,\nonumber\\
u^{\prime}&=\frac{1}{2}uv-C^2\frac{\sin\theta}{\cos^3\theta}r+r^2\,\frac{\partial V(\theta)}{\partial \theta} +\frac{\partial W(\theta)}{\partial \theta},\nonumber
\end{align}
%
%
where
\begin{align} \label{V_de_theta}
V(\theta)&=\left(\frac{M}{2}\right)^{1/2}\left(\frac{A}{\cos\theta}+\frac{4A_1}{(\cos^2\theta+\mu\sin^2\theta)^{1/2}}\right),\\
W(\theta)&=\left(\frac{M}{2}\right)^{3/2}\left(\frac{B}{\cos^3\theta}+\frac{16B_1}{(\cos^2\theta+\mu\sin^2\theta)^{3/2}}\right).  \label{W_de_theta}
\end{align}
%
%
%
%
%
%
In the new coordinates the energy integral is given by
\begin{equation}\label{e_r}
hr^3= \frac{1}{2} \left(  u^2+v^2 \right) +\frac{C^2}{2\cos^2\theta}r- r^2 V(\theta) -W(\theta)\,.
\end{equation}

\subsection{The potential functions $V(\theta)$ and $W(\theta)$}\label{V_and_W}

Recall that our study considers  the case $M>>m$ and that we introduced $\displaystyle{ \mu:=\frac{2M+m}{m}}$. In particular, we have
\begin{equation}
\label{mu_large}
\mu= 1+ \frac{2M}{m}>>1\,.
\end{equation}
In addition, we assume $\mu$ is sufficiently large so that %
\begin{equation}
\label{mu_inters_A}
\mu > 1+ \frac{A}{4A_1}\,. 
\end{equation}
A direct calculation shows that in this case
 $V(\theta)$ has three critical points at $\theta_0=0$ and $\theta=\pm\theta_v$ where
\begin{equation}
\label{theta_co}
\cos \theta_v=  \sqrt{
\frac{\mu} {(\mu-1) + (\mu-1)^{2/3} \left(\frac{4A_1}{A}  \right)^{2/3} }
}\,.
\end{equation}
Likewise, assuming that
\begin{equation}
\label{mu_inters}
\mu > 1+ \frac{B}{16B_1}\,,
\end{equation}
it follows that the function $W(\theta)$ has three critical points at $\theta_0=0$ and $\theta=\pm\theta_w$ where
\begin{equation}
\label{theta_c}
\cos \theta_w=\sqrt{\frac{ \mu}{(\mu-1)   + (\mu-1)^{2/5} \left(\frac{16 B_1}{B}\right)^{2/5}}}\,.
\end{equation}

\begin{figure}[h]
\centerline
{\includegraphics[scale=0.4]{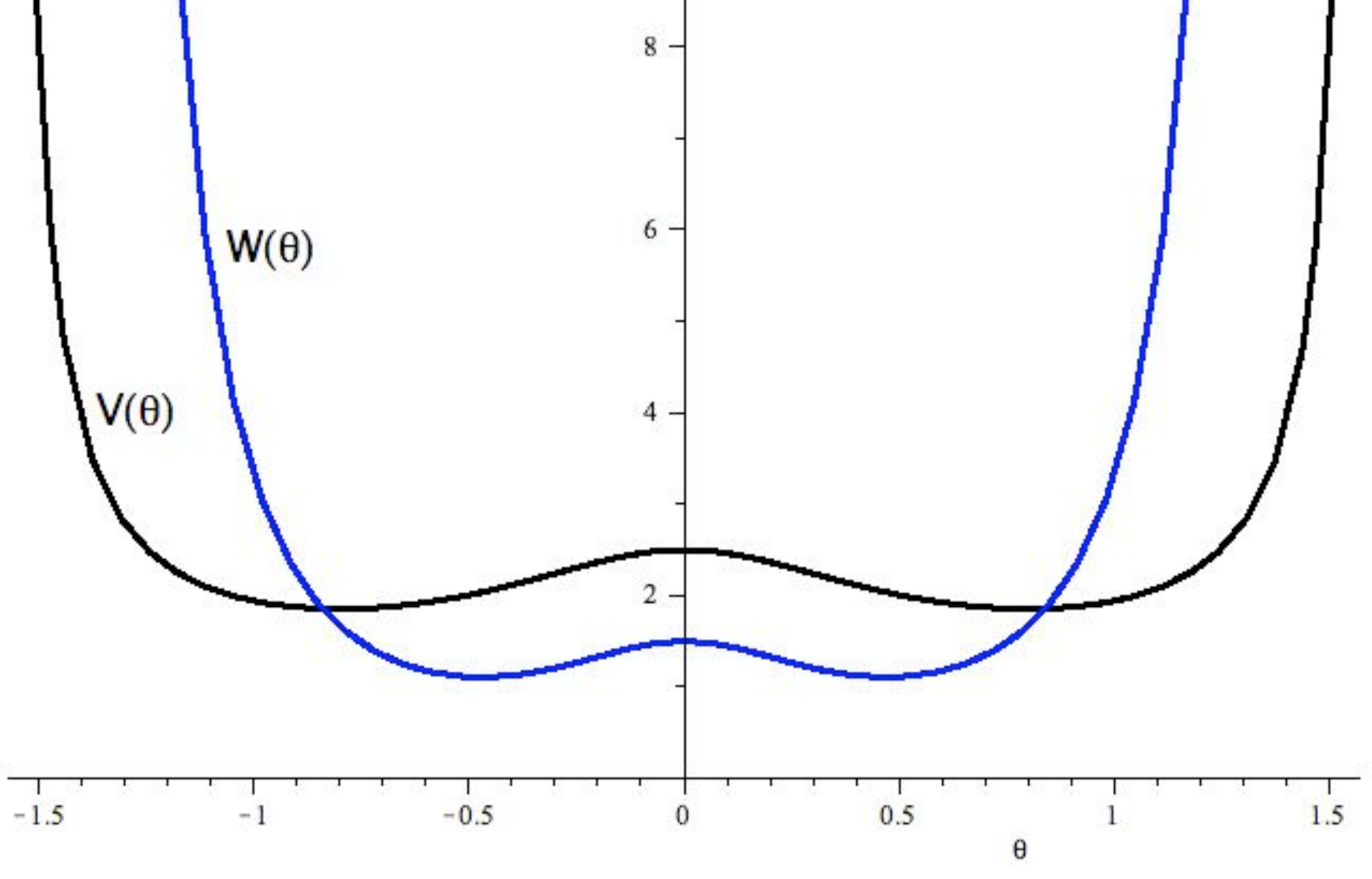}}
{\caption{ \label{shape_VW}The shape and the intersection of the functions $V(\theta)$ and $W(\theta)$ in a generic case.}}
\end{figure}

 \smallskip
 \noindent
Comparing the expressions of the non-zero critical points of $V(\theta)$ and $W(\theta),$ we deduce that in a generic situation these points do not coincide (see Figure \ref{shape_VW}). The generic case corresponds to  the condition 
\begin{equation}
\label{generic_pot}
(\mu-1)^{4/15}  \left(\frac{4A_1}{A}  \right)^{2/3}  \neq \left(\frac{16 B_1}{B}\right)^{2/5}\,.
 \end{equation}

 \noindent
 From now on, unless otherwise stated, we assume that (\ref{mu_large}), (\ref{mu_inters_A}), (\ref{mu_inters}) and (\ref{generic_pot}) are fulfilled.

%


\subsection{Regularized equations of motion and the triple collision manifold}

 The system (\ref{sys2}) is analytic for $(r, v, \theta, u) \in [0\,, \infty) \times \mathbb{R} \times \left( -\frac{\pi}{2}\,,  \frac{\pi}{2}\right)  \times \mathbb{R}$ and thus orbits at the triple collision $r=0$ are now well-defined. %
 To regularize the equations of motion at double collisions, i.e., at points with $\theta = \pm \pi/2$, we make the substitutions
\begin{equation}
\label{U-W}
U(\theta)=W(\theta)\cos^3\theta,  \quad \quad w=\dfrac{\cos^3\theta}{\sqrt{U(\theta)}}u,
\end{equation}
  and introduce a new time parametrization given by $\dfrac{\cos^3\theta}{\sqrt{U(\theta)}}=\dfrac{d\tau}{d \sigma}$.
  \begin{figure}[h]
\centerline
{\includegraphics[scale=0.6]{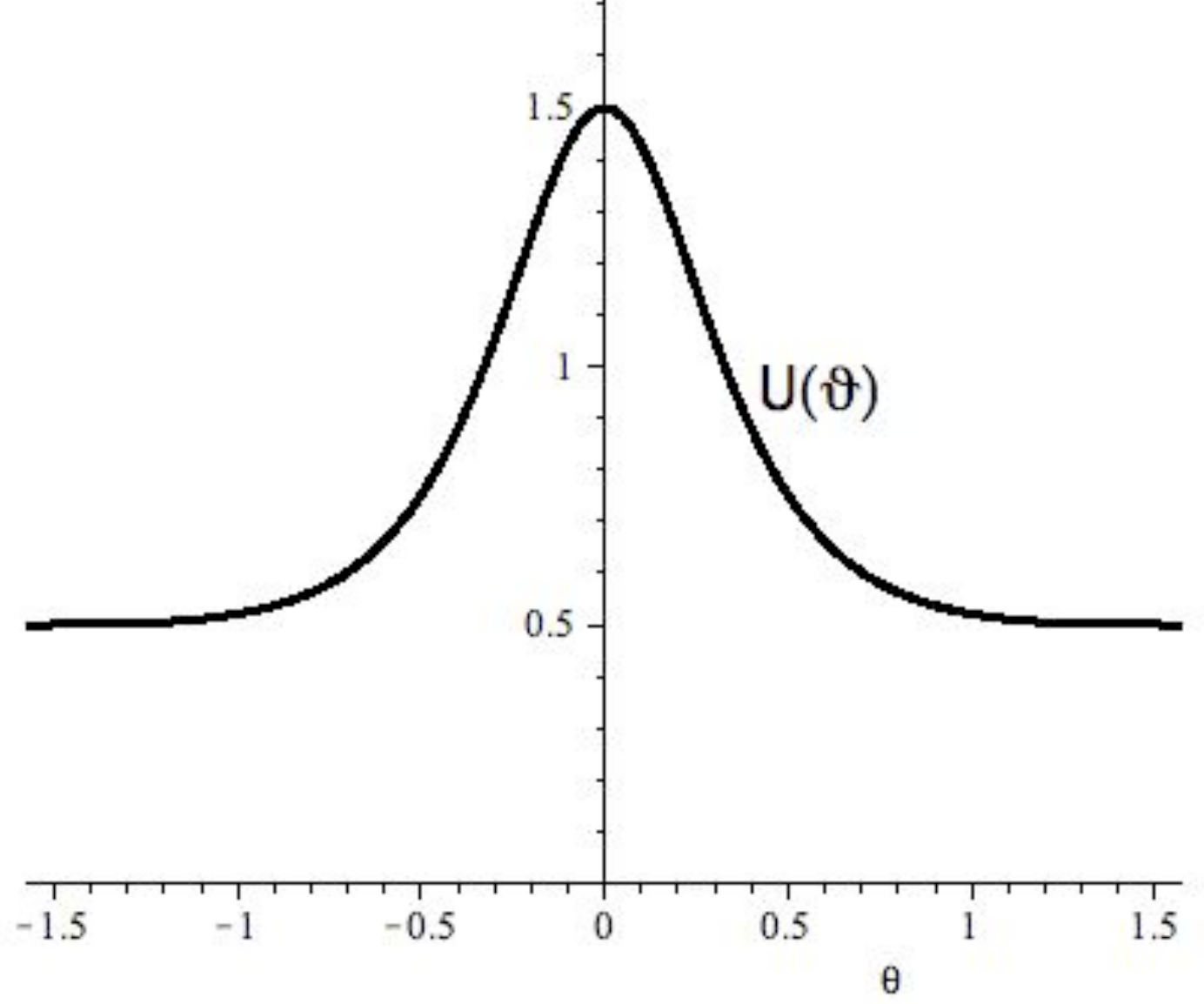}}
{\caption{ \label{shape_U}The shape of the function $U(\theta):=W(\theta)\cos^3\theta$}}
\end{figure}

\begin{figure}[h!]
\centerline
{\includegraphics[scale=0.85]{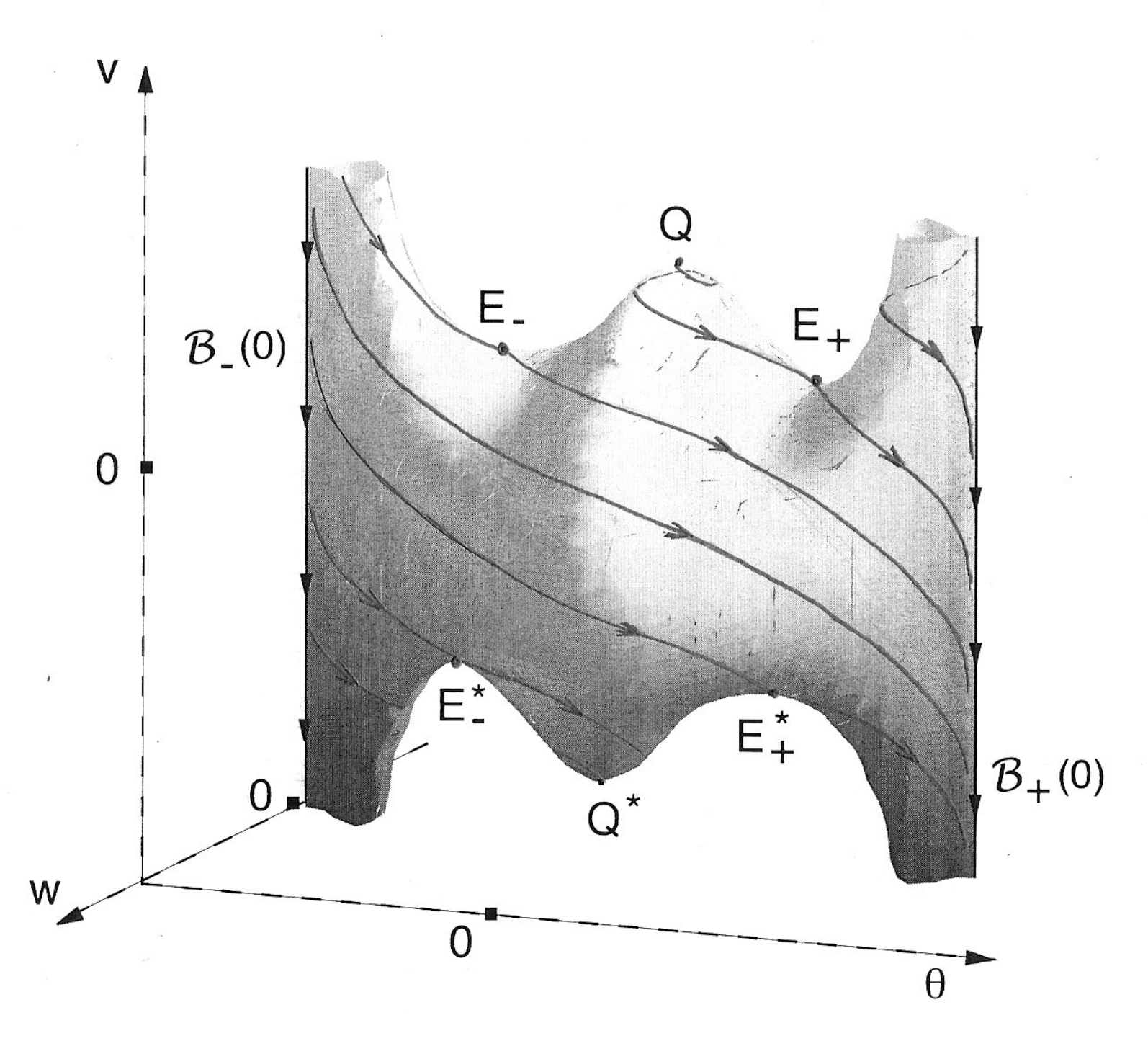}}
\caption{ \label{orbitas13}The triple collision manifold $\Delta$. The flow is gradient-like respect to the coordinate $-v$ and the sets $\mathcal{B}_{\pm}(0) $ are invariants manifolds.}
\end{figure}

\noindent
 Note that the function
 $U(\theta)>0$ for all $\theta \in [-\pi/2\,, \pi/2 ] $
and that $U(\pm\pi/2)= (M/2)^{3/2}B>0\,$ (see Figure \ref{shape_U}).
With these transformations  the system (\ref{sys2}) becomes:
\begin{align}\label{sys2_bis}
r^{\prime}&=\frac{\cos^3\theta}{\sqrt{U(\theta)}}rv,\nonumber \\
v^{\prime}&=\dfrac{\cos^3\theta}{\sqrt{U(\theta)}}\left(\frac{3}{2}v^2+\frac{U(\theta)}{\cos^6\theta}w^2-r^2V(\theta)-3\frac{U(\theta)}{\cos^3\theta} + \frac{C^2r}{\cos^2\theta} \right), \\
\theta^{\prime}&=w,\nonumber \\
w^{\prime}&=\frac{1}{2}vw\frac{\cos^3\theta}{\sqrt{U(\theta)}}+r^2V^{\prime}(\theta)\frac{\cos^6\theta}{U(\theta)}+
\frac{U^{\prime}(\theta)}{U(\theta)} \left(\cos^3\theta-\frac{w^2}{2} \right)\nonumber\\
&+3\sin\theta\cos^2\theta-\frac{\sin\theta\cos^3\theta}{U(\theta)}C^2r\nonumber,
\end{align}
where the derivation is respect to the new time $\tau$,  and the energy relation is:
\begin{equation}\label{en2}
2hr^3 \cos^6\theta= U(\theta) w^2 +\left( v^2 \cos^3 \theta - 2U(\theta) \right)\cos^3 \theta +  \left( C^2 - 2 r V(\theta)\cos^2 \theta \right) r \cos^4 \theta.
\end{equation}
Finally, using the energy relation, we substitute the term  containing the angular momentum $C$ in the $v'$ equation, and we obtain
\begin{align}\label{sys3}
r^{\prime}&=rv\, \frac{\cos^3\theta}{\sqrt{U(\theta)}},\nonumber \\
v^{\prime}&= \left(\frac{\cos^3  \theta} {2  \sqrt{U(\theta)}} v^2-\sqrt{U(\theta)} \right)+
r^2 \left(2hr +  V(\theta)\right)\frac{\cos^3\theta}{\sqrt{U(\theta)}}     , \\
\theta^{\prime}&=w,\nonumber \\
w^{\prime}&=\frac{1}{2}vw\frac{\cos^3\theta}{\sqrt{U(\theta)}}+r^2V^{\prime}(\theta)\frac{\cos^6\theta}{U(\theta)}+
\frac{U^{\prime}(\theta)}{U(\theta)} \left(\cos^3\theta-\frac{w^2}{2} \right) \nonumber\\
&+3\sin\theta\cos^2\theta -\frac{C^2r \sin\theta\cos^3\theta}{U(\theta)}\,.\nonumber
\end{align}
The vector field (\ref{sys3}) is analytic on $\displaystyle{\left[0,\infty \right)\times \R \times \left[-\frac{\pi}{2},\frac{\pi}{2} \right]\times \R}$, and thus the flow is well-defined   everywhere on its domain, including the points corresponding to triple ($r=0$) and double ($\theta=\pm \pi/2$) collisions.
The restriction of the energy relation (\ref{en2}) to $r=0$

\begin{equation}\label{en3}
\Delta:= \left\{ (r, v, \theta, w)\in \left[0,\infty \right)\times \R \times \left[-\frac{\pi}{2},\frac{\pi}{2} \right]\times \R\,|\, r=0\,, w^2+ \frac{\cos^6\theta}{U(\theta)}v^2=2\cos^3\theta \right\}
\end{equation}
 defines a fictitious  invariant manifold, called  the \textit{triple collision manifold}, pasted into the phase space for any levels of energy and angular momenta.  By continuity with respect to initial data, the flow on $\Delta$ provides  information about the orbits which pass close to collision.
The triple collision manifold is depicted in Figure \ref{orbitas13}. It is a symmetric surface with respect to the  horizontal plane $(\theta,w)$ and the  vertical plane $(v,w)$ which has on the top and bottom the profile of the function $W(\theta).$
The vector field   on the collision manifold $\Delta$ is obtained by setting $r=0$ in system (\ref{sys3}) and it is given by:
\begin{eqnarray}\label{sys4}
v^{\prime}&=&\left(\frac{\cos^3  \theta} {2  \sqrt{U(\theta)}} v^2-\sqrt{U(\theta)} \right),\nonumber \\
\theta^{\prime}&=&w,\\
w^{\prime}&=&\frac{1}{2}vw\frac{\cos^3\theta}{\sqrt{U(\theta)}}+
\frac{U^{\prime}(\theta)}{U(\theta)} \left(\cos^3\theta-\frac{w^2}{2} \right)+3\sin\theta\cos^2\theta\,.\nonumber
\end{eqnarray}

\medskip
Recall that a vector field is called \textit{gradient-like} with respect to a function $f$, if $f$ increases  along all non-equilibrium  orbits. We have:

\begin{Prop}\label{grad-like-prop}
The flow over the collision manifold is gradient-like with respect to the coordinate $-v.$
\end{Prop}

\noindent \noindent \textbf{Proof:}
On the collision manifold $\Delta$ we have
\[w^2+ \frac{\cos^6\theta}{U(\theta)}v^2=2\cos^3\theta\,.\]
Substituting $v^2$ in the expression of $v^{\prime}$ in the system (\ref{sys4}) we obtain:
\[v^{\prime}=-\frac{\sqrt{U(\theta)}}{2\cos^3\theta}w^2\]
and so $v^{\prime}<0.$ $\square$

\noindent
\begin{Prop}[Double collisions manifolds]
For each $r_0>0,$ the set 
\begin{equation} \mathcal{B}_{\pm}(r_0):=\left\{ (r, v, \theta, w) \in [0, \infty) \times  \R \times \left[-\frac{\pi}{2},\frac{\pi}{2} \right]\times \R \,|\, r=r_0\,, \theta=\pm\frac{\pi}{2}\,,w=0\right\}
\end{equation}
is an  invariant submanifold of the flow of the system (\ref{sys3}) on which the flow is gradient-like with respect to the coordinate $-v.$ \end{Prop}

\noindent \noindent \textbf{Proof:} From the equations of motion (\ref{sys3}),  for $\theta =\pm \pi/2$  we have
\begin{align}
r'&=0\,\\
v^{\prime}&=-\sqrt{U(\pm\pi/2)}<0 \,, \\
\theta^{\prime}&=w,\nonumber \\
w^{\prime}&=0\,.\nonumber
\end{align}
where we took into account that $U'(\pm\pi/2)=0\,.$   From the energy relation (\ref{en2}) we also have  that $U(\pm\pi/2)w^2=0\,,$
 from where $w \equiv 0$ and so $\theta=\pm\pi/2$ are invariant.  $\square$

\begin{Remark} 
Physically,  motions ending in $\mathcal{B}_{\pm}(r_0)$ correspond to the double collision of the masses $M$ while $m$ is located  on the vertical axis at $\displaystyle{z=\pm r_0\sqrt{\frac{2M+m}{2Mm}}}$ for $\theta=\pm \pi/2,$ respectively.
\end{Remark}

We call $\mathcal{B}_{\pm}(r_0)$  the double collision with $m$ at distance $r_0$.

\begin{Remark} Asymptotic solutions to  $\mathcal{B}_{+}(0)$ are approaching triple collision via configurations with the equal masses close to a double collision and $m$   on the same side of the vertical axis.
\end{Remark}

\medskip
 As a consequence of the previous propositions, the flow on the collision manifold consists in curves which flow  down, and either approach asymptotically the invariants sets $\mathcal{B}_{\pm}(0)$ or  end in one of the equilibrium points  (see Figure \ref{orbitas13}).



\subsection{Collision manifold  equilibria}\label{col.equi}

The equilibria of  (\ref{sys3})  at points where $r \neq 0$ are equilibria of the flow in the rotational system, and they were studied in Section \ref{Sect_RE}, Proposition\ref{real_equilibria}. 
%

\medskip
The flow on the collision manifold   accepts fictitious equilibria,  which will play an important role in understanding the orbit behavior  near singularities. We have the equilibrium points (see Figure \ref{orbitas13}):

\[Q:=(0, \sqrt{2W(0)}, 0, 0), \quad  \quad Q^*:=(0, -\sqrt{2W(0)}, 0, 0),\]
and
\[
 E_{\pm} :=(0, \sqrt{2W(\theta_w)}\,, \pm \theta_w\,, 0  ),  \quad  \quad
 E^*_{\pm} :=( 0, - \sqrt{2W(\theta_w)}\,, \pm \theta_w\,, 0  ).
\]

\bigskip
\begin{Prop}\label{dim_equi_col}
Consider the spatial Schwarzschild isosceles three body problem with parameters  such that  (\ref{mu_large}), (\ref{mu_inters_A}), (\ref{mu_inters}) and (\ref{generic_pot}) are satisfied. Let $h$ be fixed. Then on the collision manifold there are the following equilibrium points:
\[Q:=(0, \sqrt{2W(0)}, 0, 0), \quad  \quad Q^*:=(0, -\sqrt{2W(0)}, 0, 0),\]
and
\[
 E_{\pm} :=(0, \sqrt{2W(\theta_w)}\,, \pm \theta_w\,, 0  ),  \quad  \quad
 E^*_{\pm} :=( 0, - \sqrt{2W(\theta_w)}\,, \pm \theta_w\,, 0  ).
\]
Furthermore, the  equilibrium $Q$ is a spiral source with
\[ \text{dim}\, \mathcal  W_u(Q)=3,\]
 the equilibrium $Q^*$ is a spiral sink with

  \[\text{dim}\, \mathcal W_s(Q^*)=3\,. \]
The equilibria $E_{\pm}$ and $E_{\pm}^*$ are saddles with

 \[\text{dim}\,\mathcal  W_u(E_\pm)=2\,, \quad \text{dim}\,\mathcal  W_s(E_\pm)=1\,,\]

  \[\text{dim}\,\mathcal  W_s(E_\pm^*)=2\,, \quad \text{dim}\,\mathcal  W_u(E_\pm^*)=1\,.\]

\end{Prop}

\noindent
 \textbf{Proof:} The points listed above are equilibria by direct verification in the equations (\ref{sys4}). Let $\theta_c\in\{0,- \theta_w, \theta_w \}.$ The Jacobian matrix of the system (\ref{sys2_bis}) evaluating at the equilibrium point $(0, \pm \sqrt{2W(\theta_c)}, \theta_c,0)$ is
\begin{equation}J=
\left(
  \begin{array}{cccc}
  \pm  \sqrt{2\cos^3\theta_c} & 0 & 0 & 0 \\
  0 &\pm \sqrt{2}\cos^{\frac{3}{2}}\theta_c & 0 & 0 \\
    0 & 0 & 0 & 1 \\
    -\dfrac{\sin\theta_c}{W(\theta_c)}C^2 & 0 & \dfrac{W^{\prime\prime}(\theta_c)}{W(\theta_c)}\cos^3\theta_c & \pm
    \sqrt{ \frac{\cos^3\theta_c}{2} } \\
  \end{array}
\right).
\end{equation}
From the energy relation (\ref{en2}) the level of energy $h$ is given by
\begin{equation}
F(r,v,\theta,w):=-2hr^3 \cos^6\theta+ U(\theta) w^2 +\left( v^2 \cos^3 \theta - 2U(\theta) \right)\cos^3 \theta +  \left( C^2 - 2 r V(\theta)\cos^2 \theta \right) r \cos^4 \theta=0.
\end{equation}
The tangent space of this manifold at an equilibrium point $P \in \{Q, Q^*, E_{\pm}, E_{\pm}^*\}$ is
\begin{align*}
T_{P}F&=\{( \rho_1, \rho_2, \rho_3,\rho_4) \,|\, \nabla F\big|_{P} \cdot ( \rho_1, \rho_2, \rho_3,\rho_4) =0\} \\
\,\\
&= \{( \rho_1, \rho_2, \rho_3,\rho_4) \,|\, (C^2\cos^4 \theta_c)\rho_1 \pm \left(2 \sqrt{2 W(\theta_w)}\cos^6\theta_c\right)\, \rho_2=0\}.
\end{align*}

\noindent
If the angular momentum is zero, i.e., if $C=0,$ we have
$
T_{P}F=\{(\rho_1, \rho_2, \rho_3,\rho_4) \, \,|\,\rho_2=0\}.
$
The linear part of the vector field (\ref{sys2_bis}) restricted to the tangent space is given by
\begin{equation}\bar{J}=
\left(
  \begin{array}{cccc}
   \pm \sqrt{2\cos^3\theta_c}& 0 & 0 & 0 \\
    0 & 0 & 0 & 0 \\
    0 & 0 & 0 & 1 \\
    0 & 0 & \dfrac{W^{\prime\prime}(\theta_c)}{W(\theta_c)}\cos^3\theta_c &  \pm
    \sqrt{ \frac{\cos^3\theta_c}{2} } \\
  \end{array}
\right),
\end{equation}
so a basis for  $T_{P}F$ is given by  the vectors $\xi_1=(1,0,0,0)$, $\xi_3=(0,0,1,0)$ and $\xi_4=(0,0,0,1)$. A representative of $\bar{J}$ in this basis is
\begin{equation}\label{tanspace}
\left(
  \begin{array}{cccc}
   \pm \sqrt{2\cos^3\theta_c} & 0 & 0  \\
    0 & 0 & 1 \\
    0 & \dfrac{W^{\prime\prime}(\theta_c)}{W(\theta_c)}\cos^3\theta_c & \pm
    \sqrt{ \frac{\cos^3\theta_c}{2} }  \\
  \end{array}
\right).
\end{equation}
From here it follows that for  $P\in \{Q, E_{\pm}\}$ we have $\xi_1$ is a eigenvector with an eigenvalue $ \lambda_r:=\sqrt{2\cos^3\theta_c}$. For $P\in \{Q^*, E^*_{\pm}\}$ one of the eigenvalues is given by $ \lambda_r:=-\sqrt{2\cos^3\theta_c}$.   The other eigenvalues are roots of
\[\lambda^2 \mp  \sqrt{ \frac{\cos^3\theta_c}{2} } \,\lambda - \dfrac{W^{\prime\prime}(\theta_c)}{W(\theta_c)}\cos^3\theta_c =0,.  \]
with eigenvectors of the form
\[
{\bf v}^{\pm}_{\lambda_1}=
\left(
\begin{array}{c}
0\\
1\\
\lambda_1
\end{array}
  \right), \quad \quad {\bf v}^{\pm}_{\lambda_2}=
\left(
\begin{array}{c}
0\\
1\\
\lambda_2
\end{array}
  \right).
\]
The  eigenvalues at the equilibrium $Q$  $\displaystyle{\left(\text{where} \,\,\theta_c=0\,\, \text{and}\,\,v= \sqrt{2 W(0)} \right)}$ are
$$\lambda_{1,2}=
\dfrac{1}{2}\left( \frac{\sqrt{2}}{2}\pm\sqrt{\frac{25B+16B_1(1-24(\mu-1))}{2(B+16B_1)}}\right)\,,$$
%
%
and  at   $Q^*$ $\displaystyle{\left( \text{where} \,\,\theta_c=0\,\, \text{and}\,\,v=- \sqrt{2 W(0)} \right)}$  
$$\lambda_{1,2}=\dfrac{1}{2}\left( -\frac{\sqrt{2}}{2}\pm\sqrt{\frac{25B+16B_1(1-24(\mu-1))}{2(B+16B_1)}}\right)\,.$$
%
%
Given  that $\mu>> 1$, the quantity under the root is negative. It follows that  $Q$ is a spiral source and $Q^*$ is a spiral sink.
For $E_{\pm}$ $\displaystyle{\left( \text{where} \,\,\theta_c=\pm \theta_w\,\, \text{and}\,\,v= \sqrt{2 W(\theta_w)} \right)}$    the other two eigenvalues are of the form
$$\lambda_{1,2}=\frac{1}{2}\left(+ \sqrt{  \frac{\cos^3\theta_w}{2} }\pm\sqrt{\frac{\cos^3\theta_w}{2}+\frac{4W^{\prime\prime}(\theta_w)\cos^3\theta_w}{W(\theta_c)}}\right).$$
Since of $W^{\prime\prime}(\theta_w)>0$ these points are saddles. Similarly,  for $E_{\pm}^*$ 
$\displaystyle{\left( \text{where} \,\,\theta_c=\pm \theta_w\,\, \text{and}\,\,v=- \sqrt{2 W(\theta_w)} \right)}$
 we have
$$\lambda_{1,2}=\frac{1}{2}\left(- \sqrt{  \frac{\cos^3\theta_w}{2} }\pm\sqrt{\frac{\cos^3\theta_w}{2}+\frac{4W^{\prime\prime}(\theta_w)\cos^3\theta_w}{W(\theta_w)}}\right),$$
and so $E_{\pm}^*$ are saddles, too.

\bigskip
If the  angular momentum is non-zero, i.e., $C\neq 0$, then  a basis for $T_{P}F$ is given by $\xi_1=\left( \pm 2 \sqrt{2 W(\theta_w)}\cos^6\theta_c\,,C^2\cos \theta_c\,,0,0 \right)$, $\xi_3=(0,0,1,0)$ and $\xi_4=(0,0,0,1),.$
A representative of $\bar J$ in the $\{\xi_1, \xi_3, \xi_4\}$   basis is of the form
\begin{equation}\label{tanspace_2}
\left(
  \begin{array}{cccc}
     \pm \sqrt{2\cos^3\theta_c} & 0 & 0  \\
    \star & 0 & 1 \\
    \star & \dfrac{U^{\prime\prime}(\theta_c)}{U(\theta_c)}\cos^3\theta_c & \dfrac{\cos^3\theta_c}{\sqrt{U(\theta_c)}} \\
  \end{array}
\right)
\end{equation}
and the rest of the proof is identical to the one for the case $C=0$.
$\square$

\begin{cor}\label{Lebesgue}
On any energy level the set of initial conditions leading to triple collision is of positive Lebesgue measure.

\end{cor}

\begin{cor}\label{prop-RE-col-man}
For the flow restricted to  the collision manifold
  %
   the equilibrium $Q$ is a spiral source with \[ \text{dim}\,\mathcal  W_u(Q)=2,\]
 the equilibrium $Q^*$ is a spiral sink with
  \[\text{dim}\,\mathcal  W_s(Q^*)=2,\]
and the equilibria $E_{\pm}$ and $E_{\pm}^*$ are saddles with

 \[\text{dim}\,\mathcal W_s(E_\pm)=1\,, \quad \text{dim}\,\mathcal W_u(E_\pm)=1\,,\]

  \[\text{dim}\,\mathcal  W_s(E_\pm^*)=1\,, \quad \text{dim}\,\mathcal  \mathcal W_u(E_\pm^*)=1\,.\]
\end{cor}


\subsection{Orbit behavior on the triple collision manifold}\label{Orbit_triple}

On the triple collision manifold the flow is gradient like with respect to the coordinate $-v$ and has six equilibria, three in the half-space $v>0$ and three in the half-space $v<0,$ symmetrically disposed with respect to the plane $v=0.$
Given the symmetry  $\theta'(v, \theta,w) = -\theta'(v, -\theta, -w)$ and $w'(v, \theta,w) = -w'(v, -\theta, -w)$, it is sufficient to analyze the flow  on the half space $\displaystyle{\{(v, \theta, w) \in\Delta\,|\, w>0\} }$.  It is also useful to note that the vector field is invariant under $\displaystyle{\left(v(-\tau), \theta(-\tau), w(-\tau)\right) \to \left(-v(\tau), -\theta(\tau), w(\tau) \right).}$

\medskip
On the triple collision manifold the equilibrium $Q$ has a two dimensional unstable manifold. All   orbits  emerging from $Q$ flow down on $\Delta$  above $\displaystyle{\mathcal W_u(E_-)}.$
%
%
%
Looking at $E_-,$ the branch $\mathcal W_u(E_-)\big|_{\{w>0\}}$  ends either in $Q^*$, or in $E_+^*$ (and so it coincides with  $\mathcal W_s(E_+^*)\big|_{\{w>0\}}$), or  it falls in  the   basin of ${\cal B}_+(0).$  
%
%

\bigskip
In what follows we give sufficient conditions so that  $\mathcal W_u(E_-)\big|_{\{w>0\}}$   ends in  the basin of ${\cal B}_+(0)$. This will  imply that all orbits emerging from $Q$ (except for the two ending in $E_{\pm})$  flow into  ${\cal B}_+(0)$. Also, we show   that $\mathcal W_u(E_+)\big|_{\{w>0\}}$ flows into  ${\cal B}_+(0)$. This case is represented in Figure \ref{orbitas13}. 


\medskip
The flow on  the  half space $w>0$ may be obtained by substituting $w$ in the $\theta'$ equation of system (\ref{sys4}) with its expression as defined  on  the collision manifold equation (\ref{en3}). Thus, after rearranging the equation for $v'$ in (\ref{sys4}), it is given by: 
%
%
\begin{eqnarray}\label{sys7}
%
%
v^{\prime}&=&-\sqrt{U(\theta)} \left(1-\frac{\cos^3\theta}{2U(\theta)}v^2 \right),\\
\,\\
\theta^{\prime}&=&\sqrt{2\cos^3\theta\left(1 - \frac{\cos^3\theta}{2U(\theta)}v^2 \right)}. \nonumber
\end{eqnarray}
Since $\theta$ is increasing,  for $\theta \in  \left( -{\pi}/{2}\,, {\pi}/{2} \right),$   we may divide the two equations and obtain the non-autonomous differential  differential equation
\begin{equation}\label{eq_ODE}
\frac{dv}{d\theta}=
-\,\frac{1}{\sqrt{2}}\,\sqrt{ W(\theta)-\frac{1}{2}v^2 }\,.
\end{equation}
where we used that $U(\theta)= W(\theta)\cos^3 \theta$ (see (\ref{U-W})). The equation above has a smooth vector field on the domain 
\begin{equation}
\label{dom_ODE}
 {\cal D} := \left\{ (\theta, v)\,:\,   |\theta|  < \frac{\pi}{2}\,, |v| <
  \sqrt{2W(\theta)} \right\}.
\end{equation}
Also, it is symmetric under $\theta \to -\theta$ and $v \to -v.$ So whenever $v(\theta)$ is a solution, so is $-v(-\theta).$  The invariant manifold $W_u(E_-)\big|_{\{w>0\}}$ corresponds to the solution $\tilde v(\theta)$ of  (\ref{eq_ODE}) which fulfills
\begin{equation}
\label{lim_v}
\lim\limits_{\theta\to -\theta_w}\tilde v(\theta)= \sqrt{2W(-\theta_w)} =\sqrt{2W(\theta_w)} \,.
\end{equation} 
We denote by $v_1(\theta)$ the integral curve of (\ref{eq_ODE}) which passes through zero at $\theta=0$, i.e., $v_1(0)=0.$ In what follows,  we will determine  the  parameters values for which $v_1(\theta_w)>-\sqrt{2W(0)}.$ Then  
we will show that the integral curve $\tilde v(\theta)$ is above the  integral curve $ v_1(\theta)$ for all $\theta> -\theta_w$. This will imply that $\tilde v(\theta_w)>-\sqrt{2W(\theta_w)}$ and so, given that $dv/d\theta <0,$ $\tilde v(\theta)$  must tend to  $-\infty$ as $\theta \to \pi/2.$ In particular, we will obtain that  $\mathcal W_u(E_-)\big|_{\{w>0\}}={\cal B}_+(0).$ 

\medskip
We start by observing that since
 \[
 \frac{d^2v}{d\theta^2}=
 -\frac{1}{2\sqrt{2}\,\sqrt{W(\theta) -\frac{v^2}{2} } }
 \left(
 W'(\theta) - v \frac{dv}{d\theta}
 \right) =
 -\frac{1}{2\sqrt{2}\,\sqrt{ W(\theta) -\frac{v^2}{2}} }
 \left(
 W'(\theta)+\frac{v}{\sqrt{2}} \sqrt{W(\theta)- \frac{v^2}{2}}
\right)
 \]
 we have that
 \begin{equation}
 \label{concavity}
\displaystyle{  \left\{
 \begin{array}{c}
\frac{d^2v}{d\theta^2} <0 \quad \text{if} \quad  (\theta, v) \in \left\{(\theta\,, v) \,|\, \theta \in  (-\theta_w\,, 0) \,, v \in \left(0\,, \sqrt{2 W(\theta)} \right) \right\}, \\
 \,\\
 \frac{d^2v}{d\theta^2} >0 \quad \text{if} \quad  (\theta\,, v) \in \left\{(\theta, v) \,|\, \theta \in  (0\,,\theta_w) \,, v \in \left(-\sqrt{2 W(\theta)}\,, 0 \right) \right\}.
 \end{array}
 \right.}
 \end{equation}
 In other words, any integral curve of (\ref{eq_ODE}) is concave down in the upper left quadrant of $\cal D$, and concave up in the lower right quadrant of $\cal D.$

 \bigskip
 \begin{lemma}\label{lemma_ODE_1}
 In the above context, 
%
%
%
 %
if 
 \begin{equation}
 \label{cond_up}
\sqrt{\frac{ W(0)}{2} }  \leq \frac{\sqrt{2 W(\theta_w)}}{\theta_w}\,,
 \end{equation}
 then   $\tilde v(\theta)$ is well-defined for all $\theta \in ( -\theta_w, \pi/2)$ and  $\lim \limits_{\theta \to \pi/2} \tilde v(\theta)= -\infty.$ 
%

 \end{lemma}

\begin{figure}[h]
\centerline
{\includegraphics[scale=0.7]{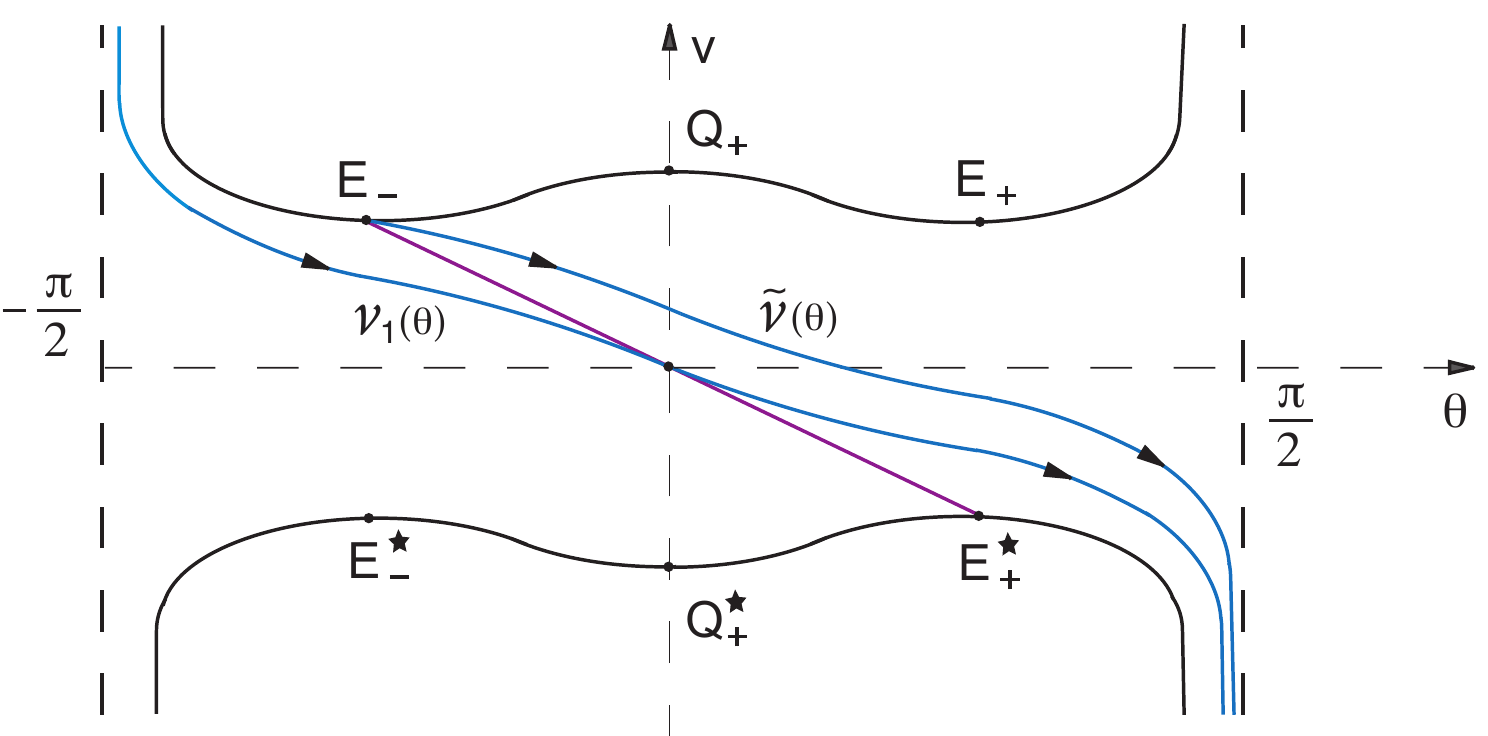}}
\caption{ \label{plano}
The domain ${\cal D}$ of the ODE (\ref{eq_ODE}) is bounded  by 
$|\theta|<\pi/2$ and $|v|<\sqrt{2W(\theta)}$. 
The solution $v_1(\theta)$ passes through $(0,0)$ where its slope is above the slope of the segment $E_-E_+^*.$ 
The solution $ \tilde v(\theta)$, which asymptotically starts in $E_-$, is always above $v_1(\theta).$
}
\end{figure}

\bigskip
\noindent
\textbf{Proof:}  Consider (\ref{eq_ODE}) and its solution which passes through $(\theta, v)=(0,0)$ which we denoted by $v_1(\theta).$ By (\ref{concavity}), $v_1(\theta)$ is concave up for $\theta >0$.
It follows that if $\displaystyle{\frac{\,\,dv_1}{d\theta} \Big|_{\theta=0}} = -\sqrt{\frac{W(0)}{2}} \geq m$, where $m$ is the slope of the segment joining $E_-$ to $E_+^*,$ then there is $\varepsilon_0>0$ such that 
 $v_1(\theta_w)=  -\sqrt{2W(\theta_w)} +\varepsilon_0 .$ 
  The inequality $\displaystyle{\frac{\,\,dv_1}{d\theta} \Big|_{\theta=0}} = -\sqrt{\frac{W(0)}{2}} \geq m$ is insured by the hypothesis condition (\ref{cond_up}). Since $v_1(\theta)$ is decreasing all along, it follows that $v_1(\theta)$  tends to ${\cal B}_+(0)$, that is $\displaystyle{ \lim\limits_{\theta\to \pi/2}v_1(\theta) =  -\infty}. $
  %
%
%
%
%
Given that $\tilde v(\theta)$ is decreasing, $\displaystyle{\tilde v(\theta)> v_1(\theta)}$ for all $\theta$  (we observe that $\tilde v (\theta)$ and $v_1(\theta)$ cannot cross as a consequence of existence and uniqueness of ODE solutions)  and $\displaystyle{ \lim\limits_{\theta\to \pi/2}v_1(\theta) =  -\infty} $, we must have 
  $\displaystyle{ \lim\limits_{\theta\to \pi/2} \tilde v(\theta) =  -\infty}$.\hspace{0.3cm}$\square$

\medskip
\begin{Remark}  The set of parameters for which (\ref{cond_up}) is fulfilled is non-empty. Indeed, after some computations (sketched below), the condition (\ref{cond_up}) is equivalent to 
 \begin{equation}
 \label{cond_param}
\cos^2 \left(
\left( \frac{4}{1+\gamma} \right)^{1/2} \left(1-\frac{1}{\mu} \right)^{3/4} \left(1+ \frac{\gamma^{2/5}}{(\mu-1)^{3/5} }   \right)^{5/4}
  \right)
\leq
 \frac{1}{\left(1-\frac{1}{\mu}\right) \left(1+  \frac{\gamma^{2/5}}{(\mu-1)^{3/5}} \right) }
 \end{equation}
 where 
$\displaystyle{\gamma:= \frac{16B_1}{B} }\,.$
%
 It can be verified (at least numerically) that for a fixed $\mu >1$, there are values of $\gamma$ which fulfill (\ref{mu_inters}) and  (\ref{cond_param}). Note that  condition (\ref{cond_param})  is independent of the angular momentum and that for $\mu \to \infty$,  at the limit it becomes $\displaystyle{ \cos^2(\sqrt{1/(1+\gamma)}) \leq 1}\,.$  %

\medskip
To see that (\ref{cond_up})  is equivalent to (\ref{cond_param}), first we substitute the definition (\ref{theta_c}) for $\theta_w$ in the expression of $W(\theta_w)$, and calculate $W(\theta_w)$, as well as $W(0)\,.$  We obtain
\[
 \theta^2_w \leq \frac{4\,W(\theta_w)}{W(0)}\,  =
\frac{
4\left( (\mu-1) + (\mu-1)^{2/5} \gamma^{2/5}  \right)^{3/2} \mu^{-3/2}\left( 1+ \gamma^{2/5} (\mu-1)^{-3/5}  \right)
}{1+\gamma}\,.
\]
After some algebra, the inequality above becomes
\[\theta^2_w  \leq \frac{4}{1+\gamma} \left(1-\frac{1}{\mu} \right)^{3/2} \left(1+ \frac{\gamma^{2/5}}{(\mu-1)^{3/5} }   \right)^{5/2} \]
which (given that $\theta_w\in (0, \pi/2)$) is equivalent to
\[\cos^2\theta_w  \geq\cos^2 \sqrt{\frac{4}{1+\gamma} \left(1-\frac{1}{\mu} \right)^{3/2} \left(1+ \frac{\gamma^{2/5}}{(\mu-1)^{3/5} }   \right)^{5/2}}\,. \]
After using again (\ref{theta_c}) and some algebra, the relation above can be written as (\ref{cond_param}).

\end{Remark}

\begin{cor}\label{match_1}
If (\ref{cond_param}) is fulfilled, then on the triple collision manifold $\displaystyle{
{\cal W}_u(E_-)\big|_{w>0}={\cal {B}}_+(0).}$
\end{cor}

\begin{cor}\label{match_2}
If (\ref{cond_param}) is fulfilled, then on the triple collision manifold all orbits emerging from $Q$ end in ${\cal {B}}_{\pm}(0)$, except for two which end in in $E_{\pm}.$ 
\end{cor}


\section{Aspects of global flow}\label{GFB.sect}

In this section we discuss homothetic solutions (defined below) and use the  properties of the flow on the collision manifold in order to analyze  the orbit behavior near double and triple collision.

\subsection{Homographic solutions}\label{hom_sol}

By definition, a \textit{homographic} solution for  a $N$ mass point system is a solution along which the  geometric configuration of the mass points is similar to the  initial geometric configuration.  If the motion of the mass points is  a uniform rotation of the initial   configuration (which stays rigid all along), then solution is in fact a  relative equilibrium.  If the mass points evolve on straight lines while forming a configuration similar to the initial configuration, then the solution is called \textit{homothetic}. In general, a homographic solution is a  superposition of a dilation and a rotation of the  initial geometric configuration of the system. Denoting by  ${\bf q}:=\left( {\bf q}_1(t), {\bf q}_2(t), \ldots, {\bf q}_N(t) \right)\in \mathbb{R}^{3N}$ the trajectories of the mass points $m_i$, it  has 
${\bf q}_i(t) = \varphi(t)   \Omega(t) \,{\bf a}_i,$ $i=1,2, \ldots, N,$ for some scalar function $\varphi(t)$ with $\text{Im}\, \varphi \in \mathbb{R}\setminus \{0\}$, some path $\Omega(t)\in SO(3)$ and some fixed non-zero vector $({\bf a}_1, {\bf a}_2,...,{\bf a}_N)\in \mathbb{R}^{3N}.$ 
If the geometric configuration (up to dilations and rotations) of the system is such that at all times the position vectors are parallel to the acceleration vectors,  then it is called a \textit{central configuration}. 

\medskip
For the isosceles three body problem   in Jacobi coordinates ${\bf r}_1=(x_1,y_1,0)$ and ${\bf r}_2=(0,0,z_2)$  homographic solutions  take the form 
 \begin{equation}
 \label{homothetic_1}
 {\bf r}_1(t)=\varphi(t)\Omega(t) {\bf a}_1\quad \text{and}\quad {\bf r}_2(t)=\varphi(t){\bf a}_2
 \end{equation}
  for some  scalar function $\varphi(t)$ with $\text{Im}\, \varphi \in \mathbb{R}\setminus \{0\}$, some rotation $\Omega(t)$ about the vertical $Oz$ axis, and a  configuration given by the fixed vectors ${\bf a}_1=(a_{1x}, a_{1y}, 0)$, and  ${\bf a}_2=(0, 0, a_{2z})\,,$ where  $({\bf a}_1, {\bf a}_2) \neq ({\bf 0},{\bf 0})$.
 Denote 
\[
\Omega(t) = \left[
\begin{array}{ccc}
\cos \psi(t) & -\sin \psi(t) & 0\\
\sin \psi(t) & \,\,\,\,\,\cos \psi(t) &0\\
0&0&1
\end{array}
\right]\,.
\]
In cylindrical coordinates $(R, \phi)$, equations (\ref{homothetic_1})  are equivalent to
\begin{align*}
R(t) \cos \phi(t) &= \varphi(t)\left( a_{1x}  \cos\psi(t) - a_{1y} \sin \psi(t) \right)\,\\
R(t) \sin \phi(t) &=\varphi(t)\left( a_{1x}  \sin\psi(t) + a_{1y}\cos \psi(t) \right)\\
 z(t) &=  \varphi(t)a_{2z}
\end{align*}
%
and so
\[R^2(t) = \varphi^2(t)(a_{1x}^2 +a_{1y}^2 ) \quad \quad \text{and} \quad \quad z(t)=  \varphi(t)a_{2z}\,.\]
Passing now to the  $(r, \theta)$  coordinates of the  equations (\ref{sys2}), we obtain that  homographic solutions satisfy
\begin{equation}
\label{cons_1}
\sqrt{\frac{2}{M}}\, r(\tau)\cos \theta(\tau)= \varphi(\tau)\sqrt{a_{1x}^2 +a_{1y}^2 }\,, \quad \quad \sqrt {\frac{2M+m}{2Mm}}\,r(\tau) \sin \theta(\tau)= \varphi(\tau)a_{2z}\,,
\end{equation}
where
\begin{equation}\label{r_not_zero}
r (\tau)\neq   0 \quad \text{for any} \quad \tau\,.
\end{equation}
From the equation of motion of (\ref{sys3}),  note that  if there is a $\tau_0$ such that $r (\tau_0) =  0$ then $r(\tau)\equiv 0$ for all $\tau.$
%
%

%
%
%

\bigskip
\noindent
If there is a $ \bar \tau$  such that $ \theta(\bar \tau) =\pi/2$ or $\theta(\bar  \tau) =-\pi/2$ then, from (\ref{cons_1}), at this value we must have $\varphi (\bar \tau)=0$ which contradicts the fact that $\text{Im}\, \varphi \in \mathbb{R}\setminus \{0\}$. Thus $\theta(\tau)\neq \pm \pi/2$ for all $\tau,$  and so we can divide equations  (\ref{cons_1})  and obtain that the homographic solutions must have 
$
\tan \theta(\tau)$ constant,
 or, equivalently 
\begin{equation}
\label{tau_cond_1}
\theta(\tau)= const.=\theta(0)=:\theta_0\,.\end{equation}
Given that $r(\tau)\neq 0$ and $\theta(\tau)\neq \pm \pi/2$ for all $\tau,$  we can work on the system (\ref{sys2}) rather than on (\ref{sys3}). For reader's convenience we re-write (\ref{sys2})  bellow:

\begin{align}\label{sys2_v2}
r^{\prime}&=rv,\nonumber \\
v^{\prime}&=\frac{3}{2}v^2+u^{2}+\frac{C^2}{\cos^2\theta}r-r^2V(\theta)-3W(\theta), \\
\theta^{\prime}&=u ,\nonumber\\
u^{\prime}&=\frac{1}{2}uv-C^2\frac{\sin\theta}{\cos^3\theta}r+r^2\,\frac{\partial V(\theta)}{\partial \theta} +\frac{\partial W(\theta)}{\partial \theta},\nonumber
\end{align}
%
%
%
%
%
%
%
with the energy relation
\begin{equation}\label{e_r_v2}
hr^3= \frac{1}{2} \left(  u^2+v^2 \right) +\frac{C^2}{2\cos^2\theta}r- r^2 V(\theta) -W(\theta)\,.
\end{equation}

\bigskip
\noindent
%
%
%
So a homographic solution is a  solution of (\ref{sys2_v2}) subject to (\ref{r_not_zero}) and (\ref{tau_cond_1}).
Now, given (\ref{tau_cond_1}), using the third equation of  (\ref{sys2_v2}) we must have     $\theta'(\tau)=u(\tau)=0$ for all $\tau.$\
%

 Let the  initial condition of a homographic solution be $(r_0, v_0, \theta_0, 0)$.  Then  $\displaystyle{\left(r(\tau), v(\tau)  \right)}$ must satisfy
\begin{align}\label{cc_101}
r^{\prime}&=rv \\
v^{\prime}&=\frac{3}{2}v^2+\frac{C^2}{\cos^2\theta_0}r-r^2V(\theta_0)-3W(\theta_0), \label{cc_102}\\
0&=-C^2\frac{\sin\theta_0}{\cos^3\theta_0}r+r^2\,V'(\theta_0) +W'(\theta_0),\label{cc_11}
\end{align}
together with the energy constraint
\begin{equation}\label{en_cc_103}
hr^3= \frac{1}{2} v^2  +\frac{C^2}{2\cos^2\theta_0}r- r^2 V(\theta_0) -W(\theta_0)\,.
\end{equation}
%
%
%
%

\bigskip
\noindent
To solve the system (\ref{cc_101})-(\ref{en_cc_103}), we distinguish the following three cases: 1) $\theta_0=0,$ 2) $\theta_0= \pm \theta_v$, i.e., $\theta_0$ is a non-zero critical point of $V(\theta)$, and 3) $\theta_0 \in (-\pi/2, \pi/2)\setminus\{0, \pm\theta_v\}$.

\bigskip
\noindent
1) If $\theta_0=0$, since $V^{\prime}(0)=W'(0)=0,$ then equation (\ref{cc_11}) is identically satisfied. It remains  to solve
\begin{align}
r^{\prime}&=rv,\label{ht_1} \\
v^{\prime}&=\frac{3}{2}v^2+C^2r-r^2V(0)-3W(0),  \label{ht_2}
\end{align}
with
\begin{equation}\label{ht_en}
hr^3= \frac{1}{2} v^2  +\frac{C^2}{2}r- r^2 V(0) -W(0)\,.
\end{equation}
The solutions of this system describe planar motions which, in the original coordinates, take place on the invariant manifold (\ref{inv_plane}) of planar motions.   We analyse this case in detail in the next section.



\bigskip
\noindent
2) If $\theta_0=\theta_v$, then from (\ref{cc_11}) it follows that
%
%
%
%
\begin{equation}\label{cc_105}
-C^2\frac{\sin\theta_v}{\cos^3\theta_v}r+W'(\theta_v) =0\,.
\end{equation}

\medskip
(a) If $C=0$ then (\ref{cc_105}) becomes
\[
W'(\theta_v) = 0\,.
\]
and so (\ref{cc_105}) is satisfied only if $\theta_v$ is a critical point of $W$ as well. This is  a  non-generic  situation which, as mentioned in Subsection \ref{V_and_W}, is not considered here.

\bigskip
(b) If $C\neq0$ then  from  (\ref{cc_105})  it follows that $r(\tau)= const.=: r_0$ for all $\tau,$ and so $r'(\tau)\equiv 0\,.$  Using (\ref{cc_101}) it follows that $v\equiv 0.$ Further, using (\ref{cc_102}), we must have
\begin{equation}\label{en_bah}
r^2_0V(\theta_v)- \left(\frac{C^2}{\cos^2\theta_v} \right)r_0+3W(\theta_v) =0\,.
\end{equation}
Since $\theta_v \neq 0$ and $r \equiv r_0$, from (\ref{cc_105}) we must have
\[C^2 =\frac{W'(\theta_v) \cos^3\theta_v }{r_0 \,\sin \theta_v }. \]
Substituting $C^2$ as above into (\ref{en_bah}), after some calculations we obtain
\begin{equation}\label{cc_bah}
r^2_0= \frac{W(\theta_v)}{V(\theta_v)}\,  \frac{\cos \theta_v}{\sin \theta_v}
\left( \frac{W'(\theta_v)}{W(\theta_v)} - 3 \frac{\sin \theta_v}{\cos \theta_v}  \right)\,.
\end{equation}
A tedious but straight-forward calculation shows that $\displaystyle{\frac{W'(\theta_v)}{W(\theta_v)} - 3 \frac{\sin \theta_v}{\cos \theta_v} <0 }$ for $\mu>1\,.$ Since all of the other terms on the right hand side are strictly positive, we obtain that $r_0^2<0$, which is a contradiction.

\bigskip
\noindent
3) $\theta_0 \in (-\pi/2, \pi/2)\setminus\{0, \pm\theta_v\}$ In this case, using  (\ref{cc_11}) 
%
%
we deduce that  $r$ is  constant. Then, since $r' \equiv 0$, and so $v\equiv 0,$  equation (\ref{cc_102}) can be written as 
\begin{equation}\label{0_2}
r^2V(\theta_0) - \frac{C^2}{\cos^2\theta_0}r+3W(\theta_0)=0. 
\end{equation}
A necessary and sufficient condition for the system (\ref{cc_11})-(\ref{0_2}) to have solutions is that  the coefficients of $r^2$ and $r$ and the free term in the two equations coincide, that is:
\[
\frac{V'(\theta_0)}{V(\theta_0)}= \tan\theta_0 =\frac{W'(\theta_0)}{3W(\theta_0)}.
\]
The first equality leads to the equation
\[ \mu \cos\theta_0 (\mu-\cos^2 \theta_0)=0\]
which, since $\mu >1$,  has no solutions.

\bigskip
In conclusion the only homographic solutions in the Schwarzschild isosceles problem are described by the solutions of the system (\ref{ht_1}) -(\ref{ht_en}). This is the subject of the next section.


\subsection{Planar motions}\label{coll-eject}

For the isosceles problem,  due to the symmetry, planar motions are homographic solutions. In our initial setting, these planar motions takes place on the plane $z\equiv 0$ and are mentioned in Remark \ref{planar_motions}. In $(r,v, \theta, u)$ coordinates, planar motions  take place on the invariant manifold
\begin{equation}
\mathcal{P}: =\left\{ (r, v, \theta, u)\,|\, \theta=0\,, u=0 \right\}
\end{equation}
and are described  given by the system (\ref{ht_1})-(\ref{ht_en}).
%
%
This system   is a 
 one degree of freedom Hamiltonian system and thus it is possible to do a full qualitative analysis of the phase space. The relative equilibria calculated in Section \ref{Sect_RE} emerge as  equilibria of  (\ref{ht_1})-(\ref{ht_2}) (scaled by a positive factor). Indeed, due to the attractive nature of the forces, all relative equilibria belong to the plane $\{z=0\}\,,$ i.e., to the invariant manifold $\mathcal{P}$. A direct calculation shows that $C_0$, the critical value of the angular momentum found in Proposition  \ref{real_equilibria}, can be written as
 \begin{equation}
 \label{C_not}
 C_0= \sqrt[4]{12 V(0)W(0)}\,,
 \end{equation}
and we have
\begin{enumerate}

 \item For $C<C_0$, the system (\ref{ht_1})-(\ref{ht_2}) has no equilibria with $r \neq 0$.
  
 \item For $C=C_0$, the system (\ref{ht_1})-(\ref{ht_2}) has a degenerate equilibrium  
 \[
r_0:= \left( \frac{C^2}{2V(0)}\,,0  \right).
 \]

 \item For $C>C_0$, the system (\ref{ht_1})-(\ref{ht_2}) has two equilibria with $r \neq 0$ located at 
 \[
r_{1}: = \left( \frac{C^2 - \sqrt{C^4-C_0^4}}{2V(0)}\,,0  \right), \quad \quad r_{2}:= \left( \frac{C^2+ \sqrt{C^4-C_0^4}}{2V(0)}\,,0  \right)\,.
 \]
The equilibrium $r_{1}$ is a saddle, whereas $ r_{2}$ is a  centre.

\end{enumerate} 
%

%

\bigskip
Another class of equilibria is given by  $
\displaystyle{(r, v) = (0, \pm \sqrt{2W(0)})\,.}$ These equilibria are independent of the angular momentum level and mark  the intersection of the triple collision manifold $\Delta$ with $\mathcal{P}$; on the Figure \ref{orbitas13} they  correspond to the points $Q$ and $Q^*$, respectively.

\bigskip
\noindent
We are  ready to sketch the phase space of (\ref{ht_1})-(\ref{ht_2}). Using the energy integral  (\ref{ht_en}),  we have
\begin{equation}\label{ht_eq}
v=\pm\sqrt{2hr^3+2V(0)r^2-C^2r+2W(0)}
\end{equation}
and deduce the classification given below.
\begin{figure}[h]
\centerline
{\includegraphics[scale=0.4]{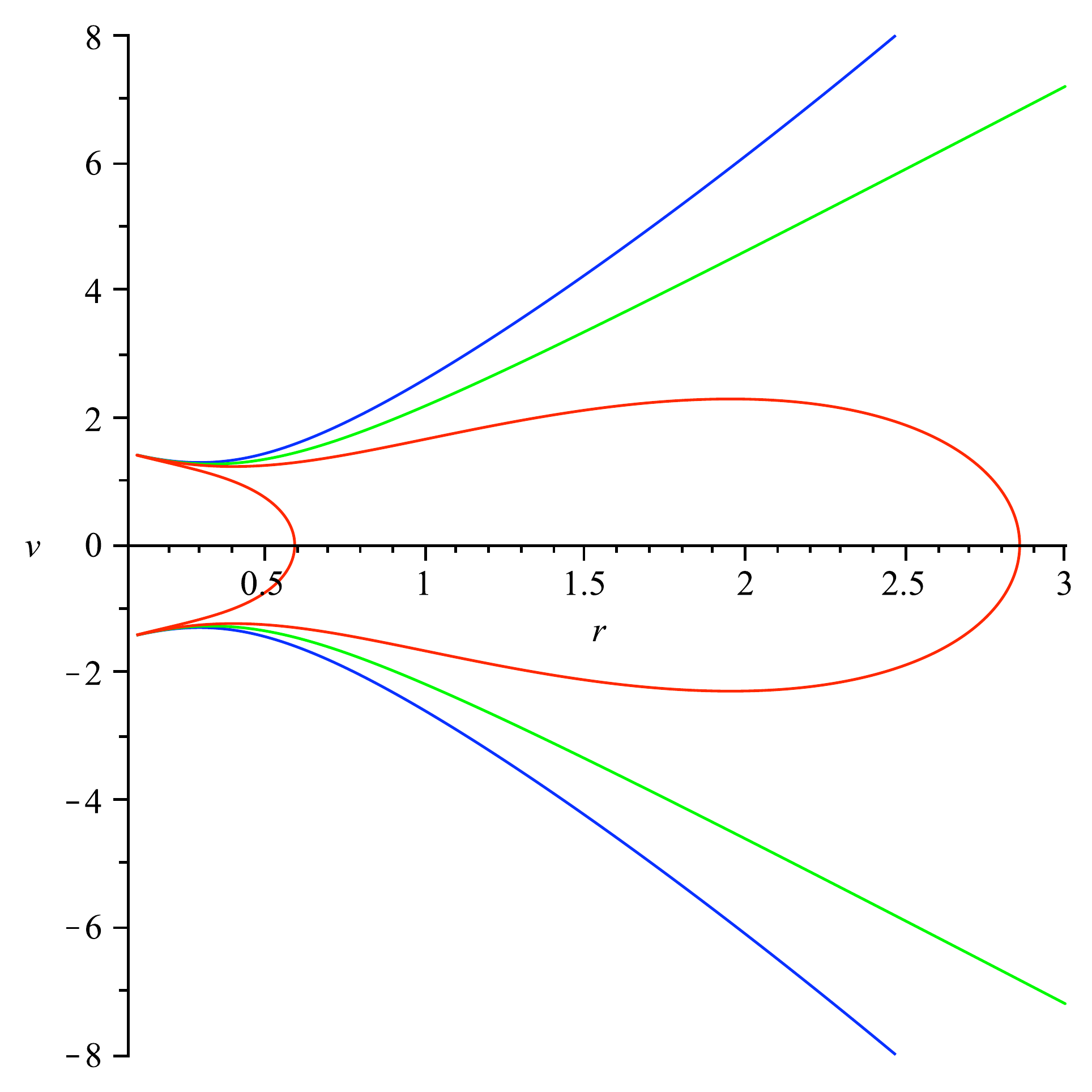}}
\caption{ \label{c_less}Curves solution of equation (\ref{ht_eq}) for a fixed  angular momentum  $0\leq C < C_0$. For $h<0$ all the orbits are bounded, for $h=0$ the orbits are parabolas and for $h>0$ cubic functions. The plot is generated for values $M=1,$ $A=A_1=1$ and $B=B_1=0.2$ and $C=C_0-0.5 \simeq 2.67-0.5\,.$  The curves sketched have (starting with the  curve from the top and going down): $h=1$ (blue),  $h=0$ (green), $h=-1$ (red) and $h=-5$ (red).
}
\end{figure}

\bigskip
\noindent
(A) For $0\leq C<C_0$, the phase curves are sketched in Figure  \ref{c_less}.  For $C=0$ all motions are rectilinear, with $m$ as the midpoints of the masses $M.$ For $0< C<C_0$ the mass points $M$ spin around the centre $m.$ We have two subcases given by $h<0$  and $h\geq0.$

\medskip
 a) For $h<0$, all orbits are bounded, ejecting from a triple collision and ending in a triple collision. For $C=0,$ the mass points $M$  eject/collide   linearly from/into $m.$ For $0< C<C_0$, the dynamics is given by  \textit{black-hole}-type orbits, where the mass points $M$ spin infinitely many times after ejecting from, and before colliding into, the third mass.  
 \begin{figure}[h]
\centerline
{\includegraphics[scale=0.4]{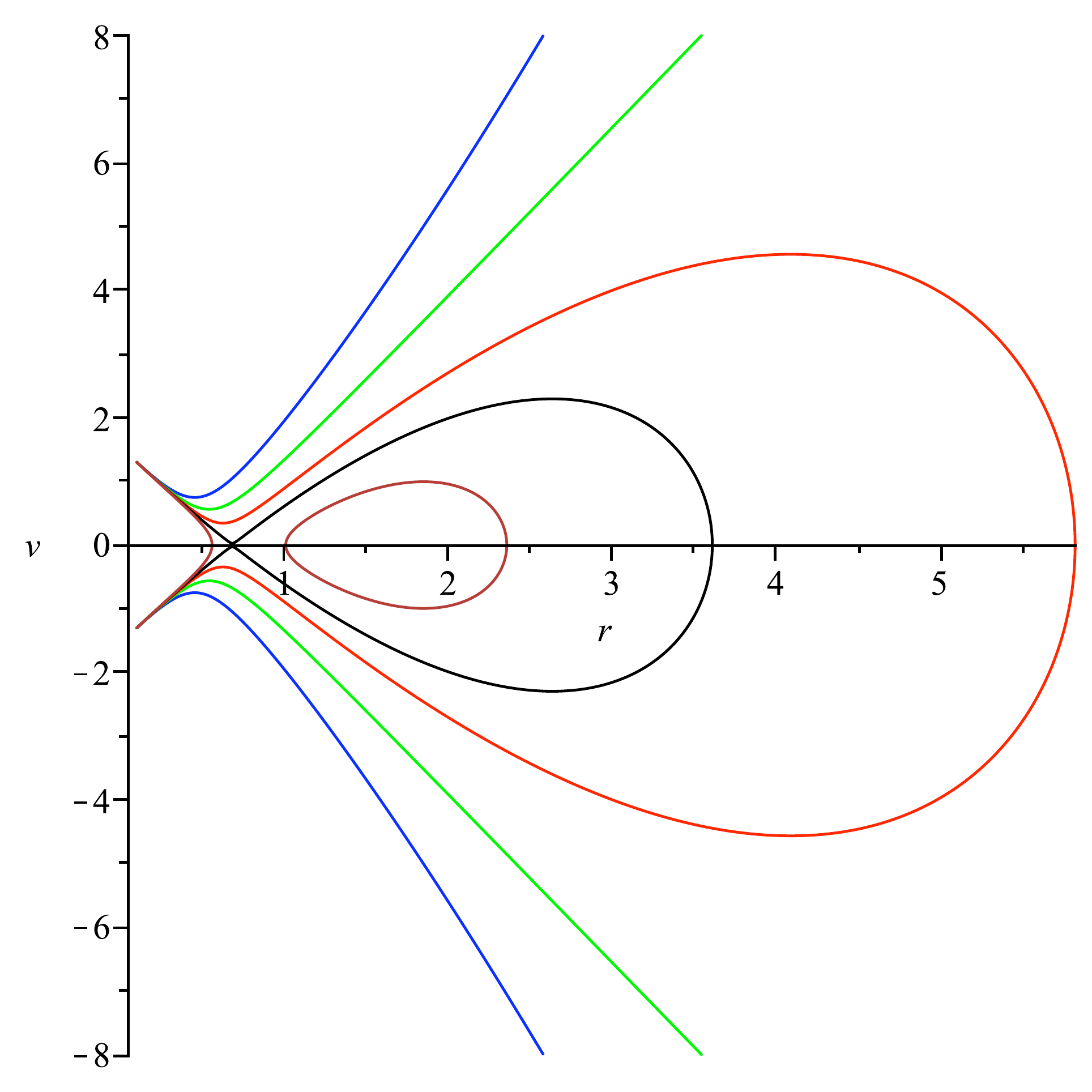}}
\caption{ \label{c_greater}Curves solution of equation (\ref{ht_eq}) for a fixed  angular momentum $C >C_0$. For $h<0$ all the orbits are bounded, for $h=0$ the orbits are parabolas and for $h>0$ cubic functions. The plot is generated for the same parameter values as in Figure \ref{c_less}.  The curves sketched have (starting with the  curve from the top and going down):   $h=1$ (blue), $h=0$ (green), $h=-0.5$ (red),  $h=-0.71$ (black), and $h=-0.9\,$ (brown) respectively.}
\end{figure} 
 
\medskip
 b) For $h\geq 0$, all orbits are unbounded, staring from triple collision and tending asymptotically to infinity. The asymptotic escape velocities at infinity  is zero and strictly positive for  $h=0$ and $h>0,$ respectively.

\bigskip
\noindent
(B) For $C=C_0$, the phase space is similar to the case $0\leq C<C_0$, except that there is a critical energy level $h=h_{cr}<0$ for which a degenerate relative equilibrium $(r_0, 0)$ lies on  the associated phase curve.

\bigskip
\noindent
(C) For $C>C_0$, we distinguish again the cases  $h<0$  and $h\geq0$ (see Figure \ref{c_greater}). 

\medskip
 a) For $h<0$  all orbits are bounded.   A large set of orbits are  ejecting from a triple collision and ending in a triple collision. Another set of orbits, bounded but of non-collisional type, is given by those surrounding the centre equilibrium $r_2$. There is also a homoclinic orbit which joins $r_1$ to itself.


  \medskip
 b) For $h\geq 0$ all orbits are unbounded. 



\section{Collision/ejection orbits}\label{n_c_o}

In this section, we describe  the orbit behavior near triple collision by employing  the information gathered on the 
    behaviour of the flow on the collision manifold  and the homographic solutions.   

\bigskip
Consider a solution of the isosceles Schwarzschild problem which asymptotically tends to a  triple collision. (An analogue reasoning can be done for a solution which asymptotically  starts in a triple collision). Such a solution must tend asymptotically  to the triple collision manifold  $\Delta$, and, in particular,  given that the flow is gradient-like on $\Delta$, to either one of the equilibria  $Q^*$, $E^*_{\pm},$ or to one of the edges $\mathcal B_{\pm}.$

Consider the collision orbit to $Q^*.$ In original coordinates, this means that 
\[\frac{z}{R} = \frac{\sqrt{\frac{2M+m}{2Mm}}  \sin \theta}{\sqrt{\frac{2}{M}}  \cos \theta} \to 0\,,\]
%
%
so the masses tend to a co-planar configuration as they approach collision. By    Section \ref{hom_sol}, any planar motion is homographic. Furthermore, the geometric configuration of such a motion is a central configuration when either the mass points are in a relative equilibrium, or the angular momentum is zero (in which case the solution is homothetic).
 Thus the limiting configuration of any triple collision orbit ending in $Q^*$
 is the geometric configuration (no necessarily central) of a homographic solution. Note that  for Newtonian interactions (see \cite{DE}) and, more generally, for  interactions given by a homogenous potential, all limiting configurations at triple collision are central configurations.
%
%

 If the triple collision orbit ends in $E^*_-$ (or $E^*_+$), then 
the limiting configuration of the mass points  forms an (non-degenerate) isosceles  triangle such that, in original coordinates,
\[\frac{z}{R} \to \tan\theta_w\,. \]
As a consequence of the analysis of Section \ref{hom_sol}, no homographic solutions have  geometric configurations of this form. So \textit{the limiting configurations of solutions which asymptotically tend to   $E^*_-$ or $E^*_+$ are  not  configurations associated to homographic solutions}. To our knowledge, this is the first time when a such ``non-homographic" configurations are observed to be  limiting configurations at triple collision.



\medskip
Given the dimensions of the stable and unstable manifolds of $Q^*$ and $E^*_{\pm}$  in Section \ref{col.equi}  we deduced that  on any energy level, the set of initial conditions leading to triple collision is of positive Lebesgue measure (see Corollary \ref{Lebesgue}).  By the remarks above, we can improve this result as follows:

\begin{Prop}
Consider the isosceles Schwarzschild problem. Then, on any energy level, the set of initial conditions leading to triple collision with a  limiting  geometric configuration of a homographic  solution is of positive Lebesgue measure. 
\end{Prop}


\begin{Prop}
Consider the isosceles Schwarzschild problem. Then, on any energy level, the set of initial conditions leading to triple collision with a non-degenerate triangular limiting configuration is of zero  Lebesgue measure.  
\end{Prop}

\medskip
 Having fixed   $h$, given that  $\text{dim}\,{\mathcal{W}}_u(Q) = \text{dim}\,{\mathcal{W}}_s(Q^*)=3$, most orbits  which pass close to    $Q$ (or $Q^*$) are in fact triple ejection-triple collision orbits which start in $Q$ and end in $Q^*$; amongst these we note the unique planar (homographic) orbit which joins $Q$ and $Q^*$.  
 This conclusion is valid for zero and non-zero angular momenta likewise. If the momentum is zero, then the motion takes place in the  vertical plane of the initial configuration, and all the  masses are just falling into  $O$,  the  midpoint of the equal masses.
 If the momentum is non-zero, the masses approach the triple collision following a scenario  where the equal masses  are on a black-hole type trajectory, spinning infinitely many times around their midpoint $O$, while  $m$ oscillates about $O$ on the vertical axis  with  decreasing amplitude.


\medskip

\medskip


Recall that the dynamics in the full space is given by the system (\ref{sys3}). We are also able to prove:

\begin{Prop}
\label{prop_global}
 Let $h<0$ be fixed. Then any solution with an initial condition   $(r_0, v_0, \theta_0, w_0)$ such that 
\begin{equation}
\label{cond_r_zero}
2 r_0^2 \tilde V(0) < \frac{\,\,C^2}{2} \quad  \text{and} \quad v_0<0\end{equation}
where
\[\tilde V(\theta) := V(\theta) \cos \theta\]
tends either to a double collision manifold ${\cal B}_{\pm}(r)$ with $0< r<r_0$, or to the triple collision manifold (including the sets ${\cal B}_{\pm}(0)$).
\end{Prop}

\noindent
\textbf{Proof:}  In  the  system (\ref{sys2_bis}), in the equation for $v'$, we substitute the term $3/2v^2$ in terms of $r,$ $\theta$ and $w$ using the energy relation (\ref{en2}), and obtain
\begin{equation}
v'(\tau)= - \frac{U(\theta)}{2 \cos^3\theta}w^2 + \frac{3h \cos^3 \theta}{\sqrt{U(\theta)}}r^3 +
\frac{r}{\sqrt{U(\theta)}} \left( 2r^2 V(\theta)\cos^2 \theta - \frac{\,\,C^2}{2}   \right)\,.
\end{equation}
Note that $\cos\theta>0$ for all $\theta \in (-\pi/2, \pi/2)$ and since $U(\theta)>0$ and $h<0,$ the coefficients of  $w^2$ and $r^3$ are negative at all times.
 
 \medskip
\noindent
 The function $V(\theta) \cos^2 \theta= \tilde V(\theta) \cos\theta $ is positive and bounded with  
 \[0< V(\theta) \cos^2 \theta\leq \tilde V(0) \quad\,\, \text{for all} \, \,\,\,\,\theta \in \left(-\frac{\pi}{2}, \frac{\pi}{2} \right).\]

\bigskip

 Consider a solution  with $(r_0, v_0, \theta_0, w_0)$ such that (\ref{cond_r_zero}) is true. Then from the  first equation of (\ref{sys2_bis}):
\[r' = \frac{\cos^3 \theta}{\sqrt{U(\theta)} }\,rv\]
 $r'_0<0$ and so $r$ is decreasing. Also,  
 \begin{align*}
v'_0&= - \frac{U(\theta_0)}{2 \cos^3\theta_0}w_0^2 + \frac{3h \cos^3 \theta_0}{\sqrt{U(\theta_0)}} r_0^3+
\frac{r_0}{\sqrt{U(\theta_0)}} \left( 2r_0^2 V(\theta_0)\cos^2 \theta_0 - \frac{\,\,C^2}{2}   \right) \\
& < - \frac{U(\theta_0)}{2 \cos^3\theta_0}w_0^2 + \frac{3h \cos^3 \theta_0}{\sqrt{U(\theta_0)}}r_0^3 +
\frac{r_0}{\sqrt{U(\theta_0)}} \left( 2r_0^2\tilde V(0)    - \frac{\,\,C^2}{2}   \right) <0\,.
\end{align*}
Thus $v$ is decreasing. At some time $\tau>0$ later, $v(\tau)<v_0<0$ and so, since $v$ is negative, $r$ will  decrease. We have $r(\tau)< r_0$ and
 \begin{align*}
v'&= - \frac{U(\theta)}{2 \cos^3\theta}w^2 + \frac{3h \cos^3 \theta}{\sqrt{U(\theta)}}r^3 +
\frac{r}{\sqrt{U(\theta)}} \left( 2r^2 V(\theta)\cos^2 \theta - \frac{\,\,C^2}{2}   \right) \\
& < - \frac{U(\theta)}{2 \cos^3\theta}w^2 + \frac{3h \cos^3 \theta}{\sqrt{U(\theta_0)}}r^3 +
\frac{r}{\sqrt{U(\theta)}} \left( 2r_0^2 V(\theta)\cos^2 \theta  - \frac{\,\,C^2}{2}   \right) 
\\
& < - \frac{U(\theta)}{2 \cos^3\theta}w^2 + \frac{3h \cos^3 \theta}{\sqrt{U(\theta_0)}} r^3+
\frac{r}{\sqrt{U(\theta)}} \left( 2r_0^2 \tilde V(0)   - \frac{\,\,C^2}{2}   \right) <0\,.
\end{align*}
 
So at $\tau$ we have fulfilled again the conditions  (\ref{cond_r_zero}) and the solution will continue to have $r$ decreasing and $v'$ negative for all $\tau.$ It follows that the solution must tend either to one of the ${\cal B}_{\pm}(r)$ for some  $r$ fixed, $r<r_0,$ or to the triple collision manifold.  $\square$

\begin{Remark}
In terms of the parameters, the value of $\tilde V(0)$ is given by

\[\tilde V(0) = \left( \frac{M}{2} \right)^{1/2} A\,.\]

\end{Remark}

\begin{Remark}
If a solution tends to a  ${\cal B}_{\pm}(r)$ with $r>0$ then, starting with a  $\tau$ large enough,  the  mass $m$ must remain on the positive or negative  side of the   vertical axis as $\theta \to \pi/2$ or $\theta \to -\pi/2$, respectively.  In particular, $m$ crosses the horizontal plane a finite number of times.

\end{Remark}

\begin{cor}
 Let $h<0$ be fixed. Then any solution with an initial condition   $(r_0, v_0, \theta_0, w_0)$ such that 

\begin{equation}
\label{cond_r_v_zero}
2 r_0^2 \tilde V(0) < \frac{\,\,C^2}{2} \quad  \text{and} \quad v_0<-\sqrt{2W(0)}\end{equation}
tends to one of the ${\cal B}_{\pm}(r)$ with $0\leq r< r_0$. If the limit is ${\cal B}_{\pm}(0),$ then the triple collision is reached after $m$ crosses the horizontal plane a finite number of times.

\end{cor}

\noindent
\textbf{Proof:}
Since $v_0< -\sqrt{2W(0)}$ and  $v$ is decreasing, the motion ends in one of the ${\cal B}_{\pm}(r)$ with $0\leq r< r_0$.  Further, note that  
for $r$ small, the only oscillatory motions with $m$ crossing the horizontal plane (or, equivalently, with $\theta$ changing its sign) are near $Q$ and $Q^*$, where $v$ is near $\pm \sqrt{2W(0)}.$ But $v_0<-\sqrt{2W(0)}$ and $v$ is decreasing. In particular,  the only possible motion tending to the collision manifold subset ${\cal B}_{\pm}(0)$ with $v$ much smaller then $- \sqrt{2W(0)}$ must have either $\theta\to \pi/2$ or $\theta\to -\pi/2$ and so, for $\tau$ large enough, oscillations are not possible anymore.

\begin{Remark}
In the conditions of Proposition \ref{prop_global}, if $C\neq0$ then the equal masses  collide following a black-hole type trajectory. 
\end{Remark}


\subsection*{Acknowledgments}

We thank the anonymous reviewers for their  useful and interesting comments and suggestions.

\medskip
The first two authors have been supported by the CONACYT Mexico grant 128790 and the third one by an NSERC Discovery Grant.

\end{document}